\documentclass[11pt, a4paper,twoside]{article}
\usepackage{amssymb,amsfonts,amsmath,amsthm,euscript}
\usepackage[english]{babel}
\usepackage[latin1]{inputenc}
\usepackage{dsfont}
\usepackage{geometry}
\usepackage{graphicx}
\usepackage{mathrsfs}
\usepackage{multirow}
\usepackage{multicol}
\usepackage[bf,footnotesize]{caption}
\usepackage{sectsty}
\usepackage{mathtools}
\usepackage{setspace}
\usepackage{placeins}
\usepackage{subfigure}

\sectionfont{\large}
\allowdisplaybreaks 
\DeclareMathAlphabet{\mz}{OT1}{pzc}{m}{it}
\DeclareMathOperator*{\diag}{\mathrm{diag}}
\DeclareMathOperator*{\tr}{\mathrm{tr}}

\begin{document}

\theoremstyle{plain}
\newtheorem{theorem}{Theorem}[section]
\newtheorem{cor}{Corollary}[section]
\newtheorem{prop}{Proposition}[section]
\newtheorem{lemma}{Lemma}[section]

\theoremstyle{definition}
\newtheorem{rmk}{Remark}[section]

\renewcommand{\theequation}{\thesection.\arabic{equation}}
\addtolength{\parskip}{.2cm}
\newcommand{\comp}[1]{\overline{#1}^{\,\,\prime}}
\newcommand{\re}{\mathrm{Re}}
\newcommand{\pim}{\mathrm{Im}}
\newcommand{\R}{\mathds{R}}
\newcommand{\N}{\mathds{N}}
\newcommand{\E}{\mathds{E}}
\newcommand{\Z}{\mathds{Z}}
\newcommand{\C}{\mathds{C}}
\renewcommand{\P}{\mathds{P}}
\renewcommand{\l}{\lambda}
\renewcommand{\L}{\Lambda}
\renewcommand{\d}{\widehat{\bs d}}
\renewcommand{\do}{\bs d_0}
\newcommand{\G}[1]{\widehat G(#1)}
\newcommand{\od}{\overline{\bs {\mathrm{d}}}}
\newcommand{\sz}{\frac{1}{m}\sum_{j=\lfloor m\kappa\rfloor}^m}
\newcommand{\qz}{\mathcal{Q}_{\kappa}(\bs d)}
\newcommand{\eps}{\varepsilon}
\newcommand{\te}{\theta}
\newcommand{\mcd}{\widehat{\mathcal{D}}}
\newcommand{\mdk}{\widehat{\mathcal{D}}_{\kappa}}
\newcommand{\mcp}{\mathcal{P}}
\newcommand{\im}{\mathrm{i}}
\newcommand{\sk}{\sum_{k=1}^q}
\newcommand{\bs}[1]{\boldsymbol{#1}}
\renewcommand{\proof}{ \noindent \textbf{\emph{Proof: }}}
\newcommand{\fim}{\hfill{\footnotesize$\blacksquare$}\normalsize\\}
\newcommand{\id}{\mathrm{I}}
\newcommand{\lfl}{\L_j(\do)^{-1}f_n(\l_j)\comp{\L_j(\do)^{-1}}}
\newcommand{\lf}{\mathcal E_j\L_j(\do)^{-1}f_n(\l_j)\comp{\L_j(\do)^{-1}}\comp{\mathcal E_j}}
\newcommand{\ir}{\mathrm{I}_{(r)}}
\newcommand{\is}{\mathrm{I}_{(s)}}
\numberwithin{equation}{section}
\numberwithin{table}{section}
\numberwithin{figure}{section}
\pagestyle{myheadings} 
\markboth{}{G. Pumi and S.R.C. Lopes} 

\long\def\sfootnote[#1]#2{\begingroup%
\def\thefootnote{\fnsymbol{footnote}}\footnote[#1]{#2}\endgroup}

\def\bfootnote{\xdef\@thefnmark{}\@footnotetext}

\thispagestyle{empty}
{\centering
\huge{\bf  A Semiparametric Estimator for Long-Range Dependent Multivariate Processes}\vspace{.8cm}\\
\large{ {\bf Guilherme Pumi$\!\phantom{i}^{\mathrm{a,}}$\sfootnote[1]{Corresponding author.}\let\thefootnote\relax\footnote{\hskip-.3cm$\phantom{s}^\mathrm{a}$Mathematics Institute - Federal University of Rio Grande do Sul - 9500, Bento Gon\c calves Avenue - 91509-900, Porto Alegre - RS - Brazil.}\let\thefootnote\relax\footnote{E-mail addresses: guipumi@gmail.com (G. Pumi), silvia.lopes@ufrgs.br (S.R.C. Lopes).} and S\'ilvia R.C. Lopes$\!\!\phantom{s}^\mathrm{a}$ } \\
\let\thefootnote\relax\footnote{This version: 12/13/2012.}\\
}
}
\vskip.6cm

\begin{abstract}
In this paper we  propose a generalization of a class of Gaussian Semiparametric Estimators (GSE) of the fractional differencing parameter for long-range dependent multivariate time series. We generalize a known GSE-type estimator by introducing some modifications at the objective function level regarding the process' spectral density matrix estimator. We study large sample properties of the estimator without assuming Gaussianity as well as hypothesis testing. The class of models considered here satisfies simple conditions on the spectral density function, restricted to a small neighborhood of the zero frequency. This includes, but is not limited to, the class of VARFIMA models. A simulation study to assess the finite sample properties of the proposed estimator is presented and supports its competitiveness. We also present an empirical application to an exchange rate data. \vspace{.2cm}\vspace{.2cm}
\end{abstract}

\section{Introduction}

Semiparametric estimation of the fractional differencing parameter in multivariate long-range dependent time series has seen growing attention in the last few years (see, for instance, Lobato, 1999, Andersen et al., 2003 and Chiriac and Voev, 2011). The first attempt to develop the theory of semiparametric estimation in the context of univariate long-range dependent time series seems to point back to late 80's with the work of K\"unsch (1987), which proposed a local Whittle-type estimator. The idea of the estimator is to locally model the behavior of the spectral density function in long-range dependent time series by locally approximating the time-domain Gaussian likelihood near the origin. These estimators comprehend the widely applied class of Gaussian Semiparametric Estimators (GSE, for short). The asymptotic theory of the particular estimator proposed in K\"unsch (1987) was challenging and posed some real theoretical difficulties in its development given the non-linear definition of the estimator. The first asymptotic results were presented in Robinson (1995b) and later further studied in Velasco (1999), Phillips and Shimotsu (2004), among others. Several variants also emerged (Hurvich and Chen, 2000, Shimotsu and Phillips, 2005, among others).

The estimation of the fractional differencing parameter in long-range dependent time series started focusing on the fully parametric case. In the univariate case, the asymptotic theory of Whittle's estimator was fully described in the work of Fox and Taqqu (1986) and Giraitis and Surgailis (1990), while the asymptotic theory of the exact maximum likelihood estimator was established by Dahlhaus (1989) and later extended to the multivariate case by Hosoya (1997), although Sowell (1989) has studied the method in the context of VARFIMA processes before. The computational cost of the exact maximum likelihood procedure in the multivariate case is high. A relatively fast approximation is studied in Luce\~no (1996) and, more recently, in Tsay (2010), both in the context of VARFIMA processes.

Fully parametric methods present some important asymptotic properties such as efficiency, $n^{1/2}$-consistency and asymptotic normality under the correct specification of the parametric model. The main criticism to the method comes from its inconsistency under misspecification of the underlying parametric structure and from the crucial role played by Gaussianity assumptions in the asymptotic theory (but see Giraitis and Surgailis, 1990), both contestable in real life applications. In this direction, the semiparametric approach presents many advantages over the parametric one, such as less distributional requirements, robustness against short-run dependencies and more efficiency compared to the latter. Another important advantage of the semiparametric approach is that  Gaussianity is usually not assumed in the asymptotic theory.

The first rigorous treatment of a multivariate semiparametric estimator was given in Robinson (1995a). Lobato (1999) analyzes a two-step GSE based on a simple local approximation of the spectral density function in the neighborhood of the origin and derive its (Gaussian) asymptotic distribution. Shimotsu (2007) analyzes another multivariate GSE by considering a refinement of the local approximation considered in Lobato (1999) and, extending the techniques in Robinson (1995b), shows its consistency and asymptotic normality. Shimotsu (2007) also considers a ``single-step'' version of Lobato (1999)'s estimator and shows its consistency and asymptotic normality. The work of Shimotsu (2007) has been recently extended to cover non-stationary multivariate long-range dependent processes in Nielsen (2011).

In this paper we are interested in semiparametric estimation of the fractional differencing parameter on multivariate long-range dependent processes. The idea is to generalize the ``single-step'' version of Lobato (1999)'s GSE considered in Shimotsu (2007) by substituting the periodogram function, applied in defining the estimator's objective function, by an arbitrary estimator of the spectral density. Although a useful tool in spectral analysis, it is well-known that the periodogram is an ill-behaved inconsistent estimator of the spectral density. Seen as a random variable, it does not even converge to a random variable at all (cf. Grenander, 1951) being considered by some authors ``\emph{an extremely poor (if not useless) estimate of the spectral density function}'' (Priestley, 1981, p.420). A natural question is can we improve the performance of the GSE estimator by considering spectral density estimators other than the periodogram? Answering this question is the main focus of our study.

Our theoretical contribution is focused on large sample properties of the proposed estimator considering different classes of spectral density estimators. First we consider objective functions with the periodogram substituted by consistent estimators of the spectral density function. We show the consistency of the GSE obtained under the same conditions in the process as considered in Shimotsu (2007). Second, we relax the consistency condition by considering spectral density estimators satisfying certain mild moment conditions and show the consistency of the related estimator. Third we examine the asymptotic normality of the proposed estimator under a certain mild regularity condition on the spectral density estimator and under the same assumptions as in Shimotsu (2007). The limiting distribution turns out to be the same as the estimator considered in Shimotsu (2007). We also consider hypothesis testing related problems. Gaussianity is nowhere assumed in the asymptotic theory.

To exemplify the use and to assess the finite sample performance of the proposed estimator, we consider the particular cases where the smoothed periodogram and the tapered periodogram are applied as spectral density estimators. We perform a Monte Carlo simulation study based on the resulting GSE and compare it to the same estimator based on the periodogram itself. We also apply the estimators to a real data set.

The paper is organized as follows. In the next section we consider some preliminary results and definitions necessary to this work and introduce the proposed estimator. In Section 3, we study the consistency of the proposed estimator while in Section 4 we derive the estimator's asymptotic distribution. In Section 5 we present a Monte Carlo simulation study in order to assess the estimator's finite sample performance and compare it to the ``single-step'' version of Lobato's estimator considered in Shimotsu (2007). The real data application is presented in Section 6 and the conclusions in Section 7. For the presentation sake, the proofs of the results in this work are presented in Appendix A.

\section{Preliminaries}
Let $\{\bs{X}_t\}_{t=0}^\infty$ be a weakly stationary $q$-dimensional process and let $f$ denote the spectral density matrix function of $\bs{X}_{\!t}$, so that
\small\[\E\big[\big(\bs{X}_{\!t}-\E(\bs{X}_{\!t})\big)\big(\bs{X}_{\!t+h}-\E(\bs{X}_{\!t})\big)'\big]=\int_{-\pi}^\pi \mathrm{e}^{\im h \l}f(\l)d\l,\]\normalsize
for $h\in\N^\ast:=\N\!\setminus\!\{0\}$. In Lobato (1999), the author considers processes for which the spectral density matrix satisfies the following local approximation
\small\begin{equation}\label{dens_spec}
f(\l)\,\sim\diag_{i\in\{1,\cdots,q\}}\{\l^{-d_i}\}\  G_0\!\!\! \diag_{i\in\{1,\cdots,q\}}\{\l^{-d_i}\}, \quad \mbox{ as } \l\rightarrow 0^+,
\end{equation}\normalsize
where $d_i\in(-0.5,0.5)$, $i=1,\cdots, q$ and $G_0$ is a symmetric positive definite real matrix.

Although specification \eqref{dens_spec} is somewhat general, as noted in Lobato (1999), several fractionally integrated models satisfy this condition.
Each coordinate process of $\{\bs{X}_t\}_{t=0}^\infty$ exhibits long-range dependence whenever the respective parameter $d_i>0$, in the sense that the respective (unidimensional) spectral density function satisfies $f(\l)\sim K\l^{-2d_i}$, as $\l\rightarrow0^+$, for  some constant $K>0$ and $i\in \{1,\cdots, q\}$.

Let
\small\begin{equation}
\label{perio}w_n(\l)\vcentcolon=\frac{1}{\sqrt{2\pi n}}\sum_{t=1}^n\bs X_{\!t} \mathrm{e}^{\im t\l}\quad\mbox{ and }\quad I_n(\l)\vcentcolon=w_n(\l)\comp{w_n(\l)}
\end{equation}\normalsize
be the discrete Fourier transform and the periodogram of $\bs X_{\! t}$ at $\l$, respectively, where $\comp{A}$ denotes the conjugate transpose of a complex matrix $A$.

From the local form of the spectral density function \eqref{dens_spec} and the frequency domain Gaussian log-likelihood localized in the neighborhood of zero, Lobato (1999) introduces a two-step  Gaussian semiparametric estimator for the parameter  $\bs d=(d_1,\cdots,d_q)^\prime$, henceforth denoted by $\widetilde{\bs d}$. The estimator is a two-step optimization procedure based on the objective function
\small\begin{equation}\label{obj.lob}
\widetilde S(\bs d):=\log\big(\det \widetilde G(\bs d)\big) -2\sum_{k=1}^qd_k\frac{1}{m}\sum_{j=1}^m\log(\l_j),
\end{equation}\normalsize
where $\l_j:=2\pi j/n$, for $j=1,\cdots,m$, are the Fourier frequencies, $m=o(n)$, with $n$ denoting the sample size, and
\small\begin{equation}\label{LOB.G}
\widetilde G(\bs d):=\frac{1}{m}\sum_{j=1}^m\re\bigg[\diag_{i\in\{1,\cdots,q\}}\{\l_j^{-d_i}\} I_n(\l_j) \diag_{i\in\{1,\cdots,q\}}\{\l_j^{-d_i}\}\bigg],
\end{equation}\normalsize
where, as usual, $\re[z]$ denotes the real part of $z\in\mathds C$.
\noindent For $i\in\{1,\cdots,q\}$, let $d_i^{Q}$ denote the univariate quasi-maximum likelihood estimate (QMLE) given in Robinson (1995b) obtained from the $i$-th coordinate process. The first step is to obtain an initial estimate of $\bs d$ by calculating the univariate QMLE for each coordinate process. Let $\bs d^Q:=(d_1^Q,\cdots,d_q^Q)^\prime$. The final estimate is obtained by calculating
\small\begin{equation}\label{LOB.est}
\widetilde{\bs d}:=\bs d^Q-\bigg(\frac{\partial^2 \widetilde S(\bs d)}{\partial \bs d\partial {\bs d}^\prime}\bigg|_{\bs d=\bs d^Q}\bigg)^{\!\!-1}\!\bigg(\frac{\partial \widetilde S(\bs d)}{\partial \bs d}\bigg|_{\bs d=\bs d^Q}\bigg).\end{equation}\normalsize
Naturally, in this case the estimator of the matrix $G_0$ in \eqref{dens_spec} is just $\widetilde G:=\widetilde G(\widetilde{\bs d})$.

Under some mild conditions, Lobato (1999) shows that, when the spectral density function follows \eqref{dens_spec}, the estimator \eqref{LOB.est} satisfies $m^{1/2}(\widetilde{\bs d}-\bs d)\overset{d}{\longrightarrow} 2\big(G_0\odot G_0^{-1}+\mathrm I_q\big)$, as $m$ tends to infinity, where $\id_q$ is the $q\times q$ identity matrix and $\odot$ denotes the Hadamard produtct.

Notice that one can consider the same estimator based on the objective function \eqref{obj.lob} as a ``single-step'' estimator by solving the $q$-dimensional optimization problem
\begin{equation}\label{LOB.sing}
\widetilde{\bs d}:=\underset{\bs d\in\Theta}{\arg\min}\{\widetilde S(\bs d)\},
\end{equation}
where $\Theta$ is the parameter space, usually some subset of $(-0.5,0.5)^q$. Estimator \eqref{LOB.sing} was considered in details in Shimotsu (2007). Arguably, a two-step procedure like the one necessary to obtain $\widetilde{\bs d}$ in \eqref{LOB.est} is computationally faster than a direct $q$-dimensional optimization procedure as \eqref{LOB.sing}. In the late 90's, a direct multidimensional optimization procedure could be troublesome considering the computational resources available for the general public at the epoch. Nowadays, however, with the recent advances in computer sciences and the development of faster CPU's, a direct optimization procedure such as \eqref{LOB.sing} represents no difficulty in practice.

Shimotsu (2007) considered a more refined local approximation for the spectral density matrix, namely
\small\begin{equation}\label{fapprox}
f(\l_j)\sim \L_j(\bs d)G\comp{\L_j(\bs d)},\quad \mbox{where} \quad \L_j(\bs d)=\diag_{k\in\{1,\cdots,q\}}\{\L_j^{(k)}(\bs d)\}\quad\mbox{and} \quad \L_j^{(k)}(\bs d)=\l_j^{-d_k}\mathrm{e}^{\im(\pi-\l_j)d_k/2},
\end{equation}\normalsize
and studied the asymptotic behavior of \eqref{LOB.sing} under \eqref{fapprox}. Under some mild conditions, the author showed the consistency and asymptotical normality of the estimator \eqref{LOB.sing} under \eqref{fapprox} even though, in this case, the estimator is based on the misspecified model \eqref{dens_spec}.

Now let $\{\bs{X}_t\}_{t=0}^\infty$ be a weakly stationary $q$-dimensional process and let $f$ denote its spectral density matrix function satisfying \eqref{fapprox}. Let $f_n$ denote an arbitrary estimator of $f$ based on the observations $\bs X_1,\cdots, \bs X_n$. Consider the objective function
\begin{equation}\label{obj.gen}
S(\bs d):=\log\big(\det \big\{\widehat G(\bs d)\big\}\big) -2\sum_{k=1}^qd_k\frac{1}{m}\sum_{j=1}^m\log(\l_j),
\end{equation}
with
\small\begin{equation}\label{gen.G}
\widehat G(\bs d):=\frac{1}{m}\sum_{j=1}^m\re\bigg[\diag_{i\in\{1,\cdots,q\}}\{\l_j^{-d_i}\} f_n(\l_j) \diag_{i\in\{1,\cdots,q\}}\{\l_j^{-d_i}\}\bigg],
\end{equation}\normalsize
where, again, $\l_j:=2\pi j/n$, for $j=1,\cdots,m$ and $m=o(n)$. Notice that \eqref{gen.G} is just \eqref{LOB.G} with $f_n$ substituting the periodogram $I_n$.
Let us define the general estimator
\begin{equation} \label{gen.est}
\widehat{\bs d}:=\underset{\bs d\in\Theta}{\arg\min}\{S(\bs d)\},
\end{equation}
where $\Theta$ denotes the space of admissible estimates, usually a subset of $(-0.5,0.5)^q$. In this work we are interested in studying the asymptotic behavior and finite sample performance of the estimator \eqref{gen.est} as a function of $f_n$. The estimator \eqref{gen.est} is a refinement of Lobato (1999) and a generalization of the results on the estimator \eqref{LOB.sing} presented in Shimotsu (2007). Our asymptotic study is divided in two main cases. First, we consider the case where $f_n$ is an arbitrary consistent estimator of the spectral density $f$. Secondly, we consider the case where $f_n$ is an arbitrary estimator of $f$ satisfying a certain moment condition. The intersection between the two cases is not empty, as we shall discuss later on.

\section{Asymptotic Theory: Consistency}\label{atc}

Let $\{\bs{X}_{\!t}\}_{t=0}^\infty$ be a weakly stationary $q$-dimensional process and let $f=(f_{rs})_{r,s=1}^q$ be its spectral density matrix satisfying \eqref{dens_spec} for a real, symmetric and positive definite matrix $G_0=(G_0^{rs})_{r,s=1}^q$. Let $\bs d_0=(d_1^0,\cdots,d_q^0)^\prime$ be the true fractional differencing vector parameter. As usual, the sup-norm is denoted by $\|\cdot\|_{\infty}$ and to simplify the notation, we shall denote the $r$-th row and the $s$-th column of a matrix $M$ by $(M)_{r\bs\cdot}$ and $(M)_{\bs\cdot s}$, respectively. Before proceed with the asymptotic theory, let us state the necessary conditions for the consistency of the estimator.

\begin{itemize}
\item Assumption \textbf{A1}. As $\l\rightarrow0^+$,
\small\[f_{rs}(\lambda)=\mathrm{e}^{\im\pi(d_r^0-d_s^0)/2}G_0^{rs}\l^{-d_r^0-d_s^0}+o(\l^{-d_r^0-d_s^0}),\quad \mbox{ for all }\,r,s=1,\cdots,q.\]\normalsize
\item Assumption \textbf{A2}. The process $\{\bs X_t\}_{t=0}^\infty$ is a causal linear process, that is,
\small\begin{equation}\label{condema}
\bs{X}_{\!t}-\E(\bs{X}_{\!t})=\sum_{k=0}^\infty A_k\bs\eps_{t-k},\quad\mbox{with}\quad\sum_{k=0}^\infty\big\|A_k\big\|_{\infty}^2<\infty,
\end{equation}\normalsize
where the innovation process $\{\bs \eps_t\}_{t\in\Z}$ is a (not necessarily uncorrelated) square integrable martingale difference, in the sense that
\small\[\E(\bs\eps_t|\mathscr{F}_{t-1})=0\quad \mbox{ and }\quad\E(\bs\eps_t\bs\eps_t^\prime|\mathscr F_{t-1})=\id_q\quad\mbox{a.s.,}\]\normalsize
for all $t\in\Z$, where $\mathscr{F}_t$ denotes the $\sigma$-field generated by $\{\bs\eps_s, s\leq t\}$. Also assume that there exist a random variable $\xi$ and a constant $K>0$ such that $\E(\xi^2)<\infty$ and $\P\big(\|\bs\eps_t\|_\infty^2>\eta\big)\leq K\P(\xi^2>\eta)$, for all $\eta>0$.
\item Assumption \textbf{A3}. With $\{A_k\}_{k\in\N}$ given in \eqref{condema}, define the function
\begin{equation}\label{ast}
A(\l)\vcentcolon=\sum_{k=0}^\infty A_k\mathrm{e}^{\im k\l}.
\end{equation}
In a neighborhood $(0,\delta)$, $\delta>0$, of the origin, we assume that $A$ is differentiable and
\small\[\frac{\partial}{\partial\l}\big(\comp{A(\l)}\big)_{r\bs\cdot }=O\Big(\l^{-1}\Big\|\big(\comp{A(\l)}\big)_{r\bs\cdot }\Big\|_{\infty}\Big), \quad \mbox{as } \l\rightarrow0^+.\]\normalsize
\item Assumption \textbf{A4}. As $n\rightarrow\infty$,
\small\[\frac{1}{m}+\frac{m}{n}\longrightarrow 0.\]\normalsize
\end{itemize}

\begin{rmk}
Assumptions \textbf{A1}-\textbf{A4} are the multivariate extensions of those in Robinson (1995b) considered in Shimotsu (2007). Assumption \textbf{A1} describes the true spectral density matrix behavior at the origin. Replacing $\mathrm{e}^{\im\pi(d_r^0-d_s^0)/2}$ by $\mathrm{e}^{\im(\pi-\lambda)(d_r^0-d_s^0)/2}$ makes asymptotically no difference, since $\lim_{\l\rightarrow 0^+}\mathrm{e}^{\im \lambda}-1=0$. Assumption \textbf{A2} regards the behavior of the innovation process which is assumed to be a not necessarily uncorrelated square integrable martingale difference uniformly dominated in probability by a square integrable random variable. Assumption \textbf{A3} is a regularity condition often imposed in the parametric case as, for instance, in Fox and Taqqu (1986) and Giraitis and Surgailis (1990). Assumption \textbf{A4} is minimal but necessary since $m$ must go to infinity for consistency, but must do so slower than $n$ in view of Assumption \textbf{A1}.
\end{rmk}

Assumptions \textbf{A1}-\textbf{A4} are local ones and only regard the spectral density behavior at the vicinities of the origin. Outside a small neighborhood of the origin, no assumption on the spectral density is made (except, of course, for the integrability of the spectral density, implied by the weak stationarity of the process).

For $\beta\in[0,1]$, let $f_n$ be a $n^{\beta}$-consistent estimator of the spectral density function\footnote{that is, $n^{\beta}(f_n-f)\overset{\P}{\longrightarrow}0$} for all $\bs d_0\in B$ where $B$ is a closed subset of $\R^q$. Since the spectral density matrix of the process $\{\bs{X}_{\!t}\}_{t=0}^\infty$ is unbounded at the zero frequency when $d_k^0\in(0,0.5)$, for some $k\in\{1,\cdots, q\}$, there is no hope in finding a consistent estimator for it in this situation. Let
\small\begin{equation}\label{set}
\Omega_\beta\vcentcolon=\bigg[-\frac{\beta}{2},0\bigg]^q\bigcap\bigg(-\frac{1}{2},0\bigg]^q\bigcap B\subseteq\bigg(-\frac12,0\bigg]^q,
\end{equation}\normalsize
Next we consider the estimator \eqref{gen.est} with $\Theta=\Omega_\beta$, that is
\begin{equation} \label{cons.est}
\widehat{\bs d}:=\underset{\bs d\in\Omega_\beta}{\arg\min}\{S(\bs d)\}.
\end{equation}

In the next theorem we establish the consistency of \eqref{cons.est} under Assumptions \textbf{A1}-\textbf{A4} considering an $n^{\beta}$-consistent spectral density function estimator. For the sake of a better presentation, the proofs of all results in the paper are postponed to Appendix A.

        \begin{theorem}\label{consistency}
        Let $\{\bs{X}_{\!t}\}_{t=0}^\infty$ be a weakly stationary $q$-dimensional process and let $f$ be its spectral density matrix. Let $f_n$ be a $n^\beta$-consistent estimator of $f$ for all $\do\in B\subseteq\R^q$, for $\beta\in[0,1]$, and let $\widehat{\bs d}$ be as in \eqref{cons.est}.  Assume that Assumptions \textbf{A1}-\textbf{A4} hold and let $\bs d_0\in\Omega_\beta$. Then, $\widehat{\bs d}\overset{\P}{-\!\!\!\longrightarrow}\bs d_0$, as $n\rightarrow \infty$.
        \end{theorem}

We now study the problem of relaxing the $n^{\beta}$-consistency of the spectral density estimator $f_n$ assumed in Theorem \ref{consistency}. We consider the class of estimators $f_n:=(f_n^{rs})_{r,s=1}^q$ satisfying, for $1\leq u<v\leq m$,
\small\begin{equation}\label{condthm}
\max_{r,s\in\{1,\cdots,q\}}\bigg\{\sum_{j=u}^v\mathrm{e}^{\im(\l_j-\pi)(d_r^0-d_s^0)/2}\l_j^{d^0_r+d^0_s}f_n^{rs}(\l_j)-G_0^{rs}\bigg\}=o_\P(1)
\end{equation}\normalsize
for $\bs d_0\in\Theta\subseteq(-0.5,0.5)^q$. A relatively simpler condition implying \eqref{condthm} is as follows:

        \begin{lemma}\label{lemma2}
        Let $\{\bs{X}_{\!t}\}_{t=0}^\infty$ be a weakly stationary $q$-dimensional process and let $f$ be its spectral density matrix and assume that assumptions \textbf{A1}-\textbf{A4} hold. Let $f_n$ be an estimator of $f$ satisfying for all $r,s\in\{1,\cdots,q\}$ and $\bs d_0\in\Theta$,
        \small\begin{equation}\label{indir}
        \E\Big(\l_j^{d_r^0+d_s^0}\Big|f_n^{rs}(\l_j)-\big(A(\l_j)\big)_{r\bs\cdot}I_{\bs\eps}(\l_j)\big(\comp{A(\l_j)}\big)_{\bs\cdot s}\Big|\Big)\:= o(1), \quad \mbox{ as }\:n\rightarrow\infty,
        \end{equation}\normalsize
        where $A$ is given by \eqref{ast} and $I_{\bs\eps}$ denotes the periodogram function associated to $\{\bs \eps_t\}_{t \in \Z}$, that is,
        \small\[I_{\bs \eps}(\l):=w_{\bs\eps}(\l)\comp{w_{\bs\eps}(\l)},\qquad\text{where}\qquad w_{\bs\eps}(\l)\vcentcolon=\frac{1}{\sqrt{2\pi n}}\sum_{t=1}^n\bs\eps_t\mathrm{e}^{\im t\lambda}.\]\normalsize
        Then, for $1\leq u<v\leq m$,
        \small\[\max_{r,s\in\{1,\cdots,q\}}\bigg\{\sum_{j=u}^v\mathrm{e}^{\im(\l_j-\pi)(d_r^0-d_s^0)/2}\l_j^{d^0_r+d^0_s}f_n^{rs}(\l_j)-G_0^{rs}\bigg\}= \mathscr{A}_{uv} +\mathscr{B}_{uv},\]\normalsize
        where $\mathscr A_{uv}$ and $\mathscr B_{uv}$ satisfy
        \small\[\E\big(|\mathscr A_{uv}|\big)=o(v-u+1)\qquad\mbox{ and }\qquad \max_{1\leq u<v\leq m}\big\{\big|v^{-1}\mathscr{B}_{uv}\big|\big\}=o_\P(1).\]\normalsize
        Thus, \eqref{condthm} holds.
        \end{lemma}

The class of estimators satisfying \eqref{condthm} is non-empty since, for instance, both, the ordinary and the tapered periodogram satisfy \eqref{indir} in the form $O(j^{-1/2}\log(j+1))$ (see lemma 1 in Shimotsu, 2007 and Section 5.3 below). Condition \eqref{condthm} is just slightly more general than \eqref{indir}. For the periodogram, \eqref{condthm} can be seen directly as well, but it is more involved (see lemma 1 in Shimotsu, 2007). Condition \eqref{condthm} plays a crucial role in replacing the $n^{\beta}$-consistency assumed in Theorem \ref{consistency}, as we shall see later. From an asymptotic point of view,  \eqref{condthm} is a sufficient condition to prove the consistency of $\d$ under Assumptions \textbf{A1-A4}. This is the content of the next theorem.

        \begin{theorem}\label{nonbeta}
        Let $\{\bs{X}_{\!t}\}_{t=0}^\infty$ be a weakly stationary $q$-dimensional process and let $f$ be its spectral density matrix. Let $f_n$ be an estimator of $f$ satisfying \eqref{condthm}, for all $r,s\in\{1,\cdots,q\}$ and $\do\in\Theta$ and consider the estimator $\d$ given by \eqref{gen.est} based on $f_n$. Assume that assumptions \textbf{A1}-\textbf{A4} hold. Then $\widehat{\bs d}\overset{\P}{-\!\!\!\longrightarrow}\do$, as $n$ tends to infinity.
        \end{theorem}

\section{Asymptotic Distribution and Hypothesis Testing}\label{asym}

Let $\{\bs{X}_{\!t}\}_{t=0}^\infty$ be a weakly stationary $q$-dimensional process, let $f=(f_{rs})_{r,s=1}^q$ be its spectral density matrix. Suppose that $f$ satisfies \eqref{dens_spec} for a real, symmetric and positive definite matrix $G_0=(G_0^{rs})_{r,s=1}^q$. Let $\bs d_0=(d_1^0,\cdots,d_q^0)^\prime$ be the true fractional differencing vector parameter. Assume that the following assumptions are satisfied
\begin{itemize}
\item Assumption \textbf{B1}. For $\alpha\in(0,2\,]$ and $r,s\in\{1,\cdots,q\}$,
\small\[f_{rs}(\l)=\mathrm{e}^{\im(\pi-\l)(d_r^0-d_s^0)/2}\l^{-d_r^0-d_s^0}G_0^{rs}+O(\l^{-d_r^0-d_s^0+\alpha}),\,\,\text{ as }\l\rightarrow0^+.\]\normalsize
\item Assumption \textbf{B2}. Assumption A2 holds and the process $\{\bs\eps_t\}_{t\in\Z}$ has finite fourth moment.
\item Assumption \textbf{B3}. Assumption A3 holds.
\item Assumption \textbf{B4}. For any $\delta>0$,
\small\[\frac{1}{m}+\frac{m^{1+2\alpha}\log(m)^2}{n^{2\alpha}}+\frac{\log(n)}{m^\delta}\longrightarrow 0, \,\,\text{ as } n\rightarrow\infty.\]\normalsize
\item Assumption \textbf{B5}. There exists a finite real matrix $M$ such that
\small\[\L_j(\do)^{-1}A(\l_j)=M+o(1),\,\,\text{ as }\l_j\rightarrow0.\]\normalsize
\end{itemize}
\vskip.2cm
\begin{rmk}\label{remark}
Assumption \textbf{B1} is a smoothness condition on the behavior of the spectral density function near the origin. It is slightly stronger than Assumption \textbf{A1} and is often imposed in spectral analysis. Assumption  \textbf{B2} imposes that the process $\{\bs X_t\}_{t=0}^\infty$ is linear with finite fourth moment. This restriction in the innovation process is necessary since the proof of Theorem \ref{norm} depends on a CLT-type result for a martingale difference sequence defined as a quadratic form involving  $\{\bs\eps_t\}_{t\in\Z}$, which must have finite variance. Assumption  \textbf{B4} is the same as assumption 4' in Shimotsu (2007). In particular, it implies that $(m/n)^b=o\big(m^{-\frac{b}{2\alpha}}\, \log(m)^{-\frac{b}{\alpha}}\big)$, for $b\neq 0$.  Assumption  \textbf{B5} is the same as assumption 5' in Shimotsu (2007) and it is a mild regularity condition in the degree of approximation of $A(\l_j)$ by $\L_j(\do)$.  It is satisfied by general VARFIMA processes.
\end{rmk}

        \begin{lemma}\label{lemma1b2}
         Let $\{\bs{X}_{\!t}\}_{t=0}^\infty$ be a weakly stationary $q$-dimensional process. Let $f$ be its spectral density matrix and assume that assumptions \textbf{B1}-\textbf{B5} hold, with \textbf{B4} holding for $\alpha=1$. Let $f_n$ be an estimator of $f$ satisfying
         \small\begin{equation}\label{cond_an}
        \max_{1\leq v\leq m}\bigg\{\sum_{j=1}^v\Big[f_n^{rs}(\l_j)-\big(A(\l_j)\big)_{r\bs\cdot}I_{\bs\eps}(\l_j)\big(\comp{A(\l_j)}\big)_{\bs\cdot s}\Big]\bigg\}= o_\P\bigg(\frac{m^2}{n^{1+|d_r^0+d_s^0|}}\bigg),
         \end{equation}\normalsize
         for all $r,s\in\{1,\cdots,q\}$ and $\do\in\Theta$. Then,
         \begin{itemize}
        \item[\rm {(a)}]
        \small\begin{equation}\label{lemma1b2a}
        \hspace{-.7cm}\max_{1\leq v\leq m}\bigg\{\max_{r,s}\!\bigg\{\!\sum_{j=1}^v\l_j^{d_r^0+d_s^0}\Big[f_n^{rs}(\l_j)- \big(A(\l_j)\big)_{r\bs\cdot}I_{\bs\eps}(\l_j)\big(\comp{A(\l_j)}\big)_{\bs\cdot s}\Big]\!\bigg\}\bigg\}=o_\P\bigg(\frac{\sqrt m}{\log(m)}\bigg);
        \end{equation}\normalsize
        \item[\rm {(b)}]
        \small\begin{equation}\label{lemma1b2b}
        \max_{1\leq v\leq m}\bigg\{\max_{r,s}\!\bigg\{\sum_{j=1}^v\l_j^{d_r^0+d_s^0}f_n^{rs}(\l_j)-\mathrm{e}^{\im(\pi-\l_j)(d_r^0-d_s^0)/2}G_0^{rs}\bigg\}\bigg\} = O_\P\bigg(\frac{m^{\alpha+1}}{n^\alpha}+\sqrt{m}\log(m)\bigg).
        \end{equation}
        \normalsize
         \end{itemize}
        \end{lemma}


The next theorem establishes the asymptotic normality of estimator \eqref{gen.est} under Assumptions \textbf{B1}-\textbf{B5} considering the class of estimators $f_n$ satisfying \eqref{cond_an}. To make the presentation simpler, let us define the matrices
\begin{equation}\label{mats}
\mathcal{E}_0:=\diag_{k\in\{1,\cdots,q\}}\Big\{\mathrm e^{\im\pi d_k^0/2}\Big\}, \quad\mathcal G_0:=\re\big[\mathcal{E}_0G_0\comp{\mathcal E_0}\big]\ \ \mbox{ and }\ \ \quad g_0:=\pim\big[\mathcal{E}_0G_0\comp{\mathcal E_0}\big].
\end{equation}

        \begin{theorem}\label{norm}
        Let $\{\bs{X}_{\!t}\}_{t=0}^\infty$ be a weakly stationary $q$-dimensional process, let $f$ be its spectral density matrix and assume that assumptions \textbf{B1}-\textbf{B5} hold, with \textbf{B4} holding for $\alpha=1$. Let $f_n$ be an estimator of $f$ satisfying \eqref{cond_an},
        for all $r,s\in\{1,\cdots,q\}$ and $\bs d_0\in\Theta$. If $\d\overset{\P}{\longrightarrow}\do$, for $\do\in\Theta$, then
        \small\[m^{1/2}(\d-\do)\overset{d}{-\!\!\!\longrightarrow} N(\bs 0, \Omega),\]\normalsize
        as $n$ tends to infinity, where
        \small\[\Omega\vcentcolon=\frac12\big(\mathcal G_0\odot\mathcal G_0^{-1}+\id_q\big)^{-1}\Sigma\ \big(\mathcal G_0\odot\mathcal G_0^{-1}+\id_q\big)^{-1},\]\normalsize
        with
        \small\[\Sigma:=\mathcal G_0\odot\mathcal G_0^{-1}+\id_q+(\mathcal G_0^{-1}g_0\mathcal G_0^{-1})\odot g_0-(\mathcal G_0^{-1}g_0)\odot(\mathcal G_0^{-1}g_0)^{\prime},\]\normalsize
        and $\id_q$ the $q\times q$ identity matrix. Furthermore, $\widehat G(\d)\overset{\P}{\longrightarrow}\mathcal G_0$.
        \end{theorem}

\begin{rmk}\label{rmkan}
A careful inspection on the proof of Theorem \ref{norm} reveals that it suffices that the estimator $f_n$ satisfies parts (a) and (b) of Lemma \ref{lemma1b2} in order to it hold, which are implied by \eqref{cond_an}.
\end{rmk}

The asymptotic variance-covariance matrix $\Omega$ takes a cumbersome form. A simple case occurs when $d_{1}^0=\cdots=d_q^0$ in which case $\mathcal G_0=G_0$ and $\Omega=2[G_0\odot G_0^{-1}-\mathrm I_q]$, the variance-covariance matrix of the limiting distribution of Lobato (1999)'s two-step estimator. Also, by Theorem \ref{norm}, $\widehat G(\d)$ is not a consistent estimator of $G_0$. However, since the $(r,s)$-th element of $\mathcal G_0$ is $\mathcal G_0^{rs}=\cos\big(\pi(d_r^0-d_s^0)/2\big)G_0^{rs}$, a consistent estimator of the matrix $G_0$ under the hypothesis of Theorem \ref{norm} is
\small\[\widehat{\mathrm G}:=\tau(\d)\odot\widehat G(\d), \qquad\mbox{ where }\qquad \tau(\d)_{rs}:=\frac1{\cos\big(\pi(\hat d_r-\hat d_s)/2\big)}\,.\]\normalsize
This result allows for hypothesis testing. First, let $\widehat \Omega$ denote the matrix defined in the same way as $\Omega$ is defined in Theorem \ref{norm}, but with $\widehat{\mathrm G}$ in place of $G_0$. It follows that, under the hypothesis of Theorem \ref{norm}, $\widehat \Omega\overset{\P}{\longrightarrow}\Omega$. Let $0<s\leq q$ and let $R$ be a $s\times q$ non-zero real matrix and $\bs\nu\in\R^s$. Consider testing a set of $s$ (independent) linear restrictions on $\do$ of the form
\small\[H_0:R\do=\bs\nu \quad \mbox{versus}\quad H_a:R\do\neq\bs\nu.\]\normalsize
Assuming the conditions of Theorem \ref{norm}, under $H_0$ the test statistics
\small\begin{equation}\label{test}
T:=m(R\d-\bs\nu)^\prime\big(R\widehat{\Omega}R^\prime\big)^{-1}(R\d-\bs\nu)
\end{equation}\normalsize
is asymptotically distributed as a $\chi_s^2$ distribution. As particular cases we have: testing the process for a common fractional differencing parameter, in which case $R=(\mathrm I_{q-1}\vdots\bs 0)-(\bs 0\vdots\mathrm I_{q-1})$ with dimension $q-1\times q$ and $\bs\nu=\bs 0$, where $\bs 0$ is a vector composed by $q-1$ zeroes; testing if the process is $I(0)$, in which case $R=\mathrm I_q$ and $\bs\nu$ is a vector of $q$ zeroes. Notice that in the particular case where $f_n$ is the periodogram $I_n$, and thus, $\d=\widetilde{\bs d}$, \eqref{test} is also valid by theorem 3 in Shimotsu (2007).

\section{Simulation Study}
In this section we present a Monte Carlo simulation study to assess the finite sample performance of the estimator \eqref{gen.est} and compare it to \eqref{LOB.sing}. In order to do that, we apply the tapered periodogram and the smoothed periodogram as the spectral density estimator $f_n$ on \eqref{gen.G}. Before presenting the simulation study, let us recall some facts on the the smoothed and tapered periodograms.

\subsection{The Smoothed Periodogram}
Let $\{\bs{X}_{\!t}\}_{t=0}^\infty$ be a weakly stationary $q$-dimensional process. Let $W_ n(\cdot)\vcentcolon=\big(W_n^{ij}(\cdot)\big)_{i,j=1}^q$ be an array of functions (called weight functions) and $\{\ell(k)\}_{k\in\N}$ be an increasing sequence of positive integers. For a Fourier frequency $\l_j$, we define the smoothed periodogram of $\bs{X}_{\!t}$ by
\small\begin{equation}\label{smoothp}
f_n(\l_j)\vcentcolon=\sum_{|k|\leq \ell(n)}W_n(k)\odot w_n(\l_{j+k})\comp{w_n(\l_{j+k})},
\end{equation}\normalsize
where $w_n(\cdot)$ is given by \eqref{perio} and $\odot$ denotes the Hadamard product. If, for some $j$ and $k$, $\l_{j+k}\notin[-\pi,\pi]$, we take $w_n$ as having period $2\pi$. In practice, at zero frequency, we use a slightly different estimative, namely,
\small\[f_n(0)\vcentcolon=\re\!\bigg[ W_n(0)\odot w_n(\l_1)\comp{w_n(\l_1)}+2\sum_{k=1}^{\ell(n)}W_n(k)\odot w_n(\l_{k+1})\comp{w_n(\l_{k+1})}\bigg].\]\normalsize
When the process $\{\bs{X}_{\!t}\}_{t=0}^\infty$ is long-range dependent, its spectral density present a pole at zero frequency, so that some authors take the summation in \eqref{smoothp} restricted to $k\neq -j$. In finite samples, however, $w_n(0)$ is always well defined since the process is finite with probability one. More details on the smoothed periodogram can be found, for instance, in Priestley (1981) and references therein.

The smoothed periodogram as defined in \eqref{smoothp} is a multivariate extension of the univariate smoothed periodogram. The use of the Hadamard product in the definition allows the use of different weight functions across different components of the process, accommodating, in this manner, the necessity often observed in practice of modeling different spectral density characteristics, including cross spectrum ones, with different weight functions. We refrain from discussing the different types of weight functions in the literature, since the subject is present in most textbooks. See, for instance, Priestley (1981), where a broad account of different weight functions, their properties and further references can be found.

In the asymptotic theory, we are only interested in sequences $\{\ell(k)\}_{k\in\N}$ and weight functions $W_n(\cdot)$ satisfying
\begin{itemize}
\item Assumption \textbf{C1}. $1/\ell(n)+\ell(n)/n\longrightarrow 0$, as $n$ tends to infinity;
\item Assumption \textbf{C2}. $W_n^{ij}(k)=W_n^{ij}(-k)$ and $W_n^{ij}(k)\geq 0$, for all $k$;
\item Assumption \textbf{C3}. $\sum_{|k|\leq \ell(n)}W_n^{ij}(k)=1$;
\item Assumption \textbf{C4}. $\sum_{|k|\leq \ell(n)}W_n^{ij}(k)^2\longrightarrow 0$, as $n$ tends to infinity.
\end{itemize}
Under Assumptions \textbf{C1}-\textbf{C4}, the smoothed periodogram is an $n^{1/2}$-consistent estimator of the spectral density matrix for all $\bs d\in[-0.5,0]^q$. Assumptions \textbf{C1}-\textbf{C4} are standard ones in the asymptotic theory of the smoothed periodogram (see, for instance, Priestley, 1981). Assumption \textbf{C1} controls the convergence rate of the sequence $\{\ell(k)\}_{k\in\N}$ with respect to $n$. Assumptions \textbf{C2}-\textbf{C3} impose that the weight functions must be non-negative symmetric functions and that the sequences $\big\{W_n^{ij}(k)\big\}_{k=-\ell(n)}^{\ell(n)}$ form a convex sequence of coefficients for each $n$, $i$ and $j$. Assumption \textbf{C4} is just a mild technical condition for the consistency of the estimator.

Since the smoothed periodogram (under Assumptions \textbf{C1}-\textbf{C4}) is an $n^{1/2}$-consistent estimator of the spectral density function for $\do\in B:=[-1/2,0]^q$, Theorem \ref{consistency} applies and we conclude that the estimator \eqref{cons.est} (under Assumptions \textbf{A1}-\textbf{A4}) is consistent for all $\do\in\Omega_{1/2}=[-1/4,0]^q$. To this moment, we were not able to establish the consistency of the estimator \eqref{gen.est} based on the smoothed periodogram via Theorem \ref{nonbeta} nor its asymptotic normality via Theorem \ref{norm}. However, there is empirical evidence (as we shall show later) that this is indeed the case.

\subsection{The Tapered Periodogram}
The main idea on the tapered periodogram is to obtain a decrease on the asymptotic bias by tapering the data before calculating the periodogram of the series. This is specially helpful in the case of long-range dependent time series. See, for instance, Priestley (1981) and Hurvich and Beltr\~ao (1993).

Let $\{\bs{X}_{\!t}\}_{t=0}^\infty$ be a weakly stationary $q$-dimensional process. For $i\in\{1,\cdots,q\}$, let $h_{i}:[0,1]\rightarrow\R$ be a collection of functions. Consider the vector of functions $L_n(\cdot):=\big(L_n^{i}(\cdot)\big)_{i=1}^q$ defined as $L_n^{i}(\l):=h_{i}\big(\l/n\big)$ and let
\small\[\ S_n(\l):=\Bigg(\frac{L_n^{i}(\l)}{\sqrt{\sum_{t=1}^nL_n^{i}(t)^2}}\Bigg)_{i=1}^q.\]\normalsize
The tapered periodogram $I_T(\l;n)$ of $\{\bs X_{\!t}\}_{t=1}^n$ is then defined by setting
\small\begin{equation}\label{taper}
I_T(\l;n):=w_T(\l;n)\comp{w_T(\l;n)},\quad\mbox{ where }\quad w_T(\l;n):=\frac1{\sqrt{2\pi}}\sum_{t=1}^nS_n(t)\odot \bs X_{\!t}\mathrm{e}^{-it\l} .
\end{equation}\normalsize
We shall assume the following:
\begin{itemize}
\item Assumption \textbf{D}. The tapering functions $h_{i}$ are of bounded variation and $H_{i}:=\int_0^1h_{i}^2(x)dx>0$, for all $i\in\{1,\cdots,q\}.$
\end{itemize}
The tapered periodogram is not a consistent estimator of the spectral density function, since the reduction on the bias induces, in this case, an augmentation of the variance. Just as the ordinary periodogram, the increase in the variance can be dealt by smoothing the tapered periodogram in order to obtain a consistent estimator of the spectral density function in the case $\bs d\in(-0.5,0)$ (see, for instance, the recent work of Fryzlewicz, Nason and von Sachs, 2008). Usually, a good performance of the tapered periodogram is obtained through tapering functions which decay faster than the Féjer's kernel. For more information on the choices of taper functions, see Priestley (1981), Dahlhaus (1983), Hurvich and Beltr\~ao (1993), Fryzlewicz, Nason and von Sachs (2008) and references therein.

Under Assumption \textbf D, $\sum_{t=1}^nL_n^{i}(t)^2\sim nH_{i}$ (cf. Fryzlewicz, Nason and von Sachs, 2008) so that $I_T(\l;n)=O\big(I_n(\l)\big)$. This allows to show that the estimator  \eqref{gen.est} based on the tapered periodogram is also consistent and asymptotically normally distributed. These are the contents of the next Corollaries.

        \begin{cor}\label{tapered.cons}
        Let $\{\bs{X}_{\!t}\}_{t=0}^\infty$ be a weakly stationary $q$-dimensional process with spectral density $f$ satisfying Assumptions \textbf{A1}-\textbf{A4}. Let $f_n$ be the tapered periodogram defined in \eqref{taper} satisfying Assumption \textbf{D}. For $\do\in \Theta$, consider the estimator $\widehat{\bs d}$ based on $f_n$, as given in \eqref{gen.est}. Then, $\widehat{\bs d}\overset{\P}{-\!\!\!\longrightarrow}\do$, as $n$ tends to infinity.
        \end{cor}

        \begin{cor}\label{tapered.an}
        Let $\{\bs{X}_{\!t}\}_{t=0}^\infty$ be a weakly stationary $q$-dimensional process with spectral density $f$ satisfying Assumptions \textbf{B1}-\textbf{B5}, with \textbf{B4} holding for $\alpha=1$. Let $f_n$ be the tapered periodogram given in \eqref{taper} satisfying Assumption \textbf{D}. For $\do\in \Theta$, consider the estimator $\widehat{\bs d}$ based on $f_n$, as given in \eqref{gen.est}. Then, for $\do\in\Theta$,
        \small\[m^{1/2}(\d-\do)\overset{d}{-\!\!\!\longrightarrow} N(\bs 0, \Omega),\]\normalsize
        as $n$ tends to infinity, with $\Omega$ as given in Theorem \ref{norm}.
        \end{cor}

\subsection{Simulation Results}\label{sim}
Recall that the class of VARFIMA$(p,\bs d,q)$ processes comprehend $q$-dimensional stationary processes $\{\bs X_t\}_{t \in \Z}$ which satisfy the difference equations
\small\[\bs{\Phi}({\cal B})\,\mathrm{diag}\big\{(1-{\cal B})^{\bs d}\big\}\big(\bs X_{t}-\E(\bs X_{t})\big)=\bs{\Theta}({\cal B})\bs{\varepsilon}_{t},\]\normalsize
where $\cal B $ is the backward shift operator, $\{\bs\eps_{t}\}_{t\in\Z}$ is an $m$-dimensional stationary process (the innovation process), $\bs{\Phi}({\cal B})$ and $\bs{\Theta}({\cal B})$ are $m\times m$ matrices in ${\cal B}$, given by the equations 
\small\begin{eqnarray}
\bs{\Phi}({\cal B})=\sum_{\ell=0}^{p}\bs{\phi}_{\ell}{\cal B}^\ell\,\,\,\,
\mbox{and}\,\,\,\,
\bs{\Theta}({\cal B})=\sum_{r=0}^{q}\bs{\theta}_{r}{\cal B}^r,\nonumber
\end{eqnarray}\normalsize
assumed to have no common roots, where $\bs{\phi}_{1},\cdots,\bs{\phi}_{p}, \bs{\theta}_{1},\cdots, \bs{\theta}_{q}$ are real $m\times m$ matrices and $\bs{\phi}_0 = \bs{\theta}_0 = \bs I_{m\times m}$, the $m\times m$ identity matrix.

The simulation study is based on bidimensional Gaussian VARFIMA$(0,\bs d,0)$ time series (i.i.d. innovation process) of sample size $n=1,000$ for four different pairs of the parameter $\bs d$ and (within component) correlation $\rho\in\{0,0.3,0.6,0.8\}$. A total of 1,000 replications is performed for each set of parameters. The time series are generated by the widely applied method of truncating the multidimensional infinite moving average representation of the process. The truncation point is fixed at 50,000 for all cases. The goal is the estimation of the parameter $\bs d$. To do that, we consider the estimator \eqref{gen.est} with the tapered and the smoothed periodogram as the spectral density matrix estimator $f_n$. For the tapered periodogram, we apply the cosine-bell tapering function, namely,
\small\[h_{i}(u)=\left\{\begin{array}{cc}
                \frac12\big[1-\cos(2\pi u)\big], & \mbox{ if }\ 0\leq u\leq1/2, \vspace{.2cm}\\
                h_{i}(1-u), & \mbox{ if }\ 1/2<u\leq 1,
              \end{array}\right. \quad \mbox{ for all } i=1,2. \]\normalsize
The resulting estimator is denoted by \emph{TLOB}. The cosine-bell taper is often applied as tapering function in applications as, for instance, in Hurvich and Ray (1995), Velasco (1999) and Olbermann et al. (2006). For the smoothed periodogram, we apply the so-called Bartlett's weights for all spectral density components, namely
%
\small\[W_n^{ij}(k):=\frac{\sin^2\big(\ell(n)\l_k/2\big)}{n\ell(n)\sin^2(\l_k/2)},\quad \mbox{ for all } i,j=1,2.\]\normalsize
We consider the smoothed periodogram with and without the restriction $k\neq -j$ in \eqref{smoothp} and the resulting estimator are denoted by \emph{SLOB} and $SLOB^\ast$, respectively. We also apply, for comparison purposes, the estimator given in \eqref{LOB.sing}, denoted by \emph{LOB}. The specific truncation point of the smoothed periodogram function is of the form $\ell(n):=\ell(n,\beta)=\lfloor n^{\beta}\rfloor$, for $\beta\in\{0.7,0.9\}$ while the truncation point of the objective function \eqref{obj.gen} is of the form $m:=m(n,\alpha)=\lfloor n^{\alpha}\rfloor$, for $\alpha\in\{0.65, 0.85\}$ for all estimators. All simulations were performed by using the computational resources of the (Brazilian) Center of Super Computing (CESUP-UFRGS). The routines were all implemented in FORTRAN 95 language optimized with OpenMP directives for parallel computing.

Table \ref{tab+1} and \ref{tab+2} report the simulation results for \small $\bs d\in\big\{(0.2,0.3),(0.1,0.4),(0.3,0.4),(0.1,0.3)\big\}$\normalsize.  Presented are the estimated values (mean), their standard deviations (st.d.) and the mean square error of the estimates (mse).

\begin{table}[!h]
\renewcommand{\arraystretch}{1.2}
\setlength{\tabcolsep}{3pt}
\caption{Simulation results for the estimators \emph{SLOB}, $SLOB^\ast$, \emph{TLOB} and \emph{LOB} in Gaussian VARFIMA$(0,\bs d,0)$ processes. Presented are the estimated values (mean), their standard deviation (st.d) and mean square error of the estimates (mse).}\label{tab+1}
\vskip.3cm
\centering
{\scriptsize
\begin{tabular}{|c|c|c|c|ccc|ccc|ccc|ccc|}
\hline \hline
\multirow{3}{*}{$\rho$}&\multirow{3}{*}{Method}&\multirow{3}{*}{$\beta$}&   \multirow{3}{*}{$\hat d_i$}    &\multicolumn{6}{|c|}{$\bs d=(0.1,0.4)$}
      &\multicolumn{6}{c|}{$\bs d=(0.2,0.3)$}   \\
      \cline{5-16}
      &       &       &      &\multicolumn{3}{|c|}{$\alpha=0.65$}       &\multicolumn{3}{|c|}{$\alpha=0.85$}            &\multicolumn{3}{|c|}{$\alpha=0.65$}
      &\multicolumn{3}{|c|}{$\alpha=0.85$}    \\
      \cline{5-16}
      &       &       &      & mean & st.d. & mse& mean & st.d. & mse& mean & st.d. & mse& mean & st.d. & mse\\
      \hline\hline
      \multirow{12}{*}{$0$}&\multirow{4}{*}{  \emph{SLOB}   }&\multirow{2}{*}{ 0.7 }  &$\hat d_1$&0.1024 & 0.0543 & 0.0030 &  0.0954 & 0.0268 & 0.0007 &  0.2081 & 0.0569 & 0.0033 &  0.1944 & 0.0277 & 0.0008\\
&& & $\hat d_2$ &  0.4470 & 0.0788 & 0.0084 &  0.4138 & 0.0401 & 0.0018 &  0.3167 & 0.0644 & 0.0044 &  0.2999 & 0.0321 & 0.0010\\
\cline{3-16}
      &&\multirow{2}{*}{ 0.9   }  &$\hat d_1$&0.1044 & 0.0560 & 0.0032 &  0.0953 & 0.0267 & 0.0007 &  0.2085 & 0.0564 & 0.0033 &  0.1923 & 0.0269 & 0.0008\\
&& & $\hat d_2$ &  0.4233 & 0.0624 & 0.0044 &  0.3950 & 0.0306 & 0.0010 &  0.3080 & 0.0580 & 0.0034 &  0.2919 & 0.0289 & 0.0009\\
\cline{2-16}
      &\multirow{4}{*}{  \emph{SLOB}\hspace{-.25cm}$\phantom{B}^\ast$   }&\multirow{2}{*}{ 0.7 }  &$\hat d_1$&0.0922 & 0.0549 & 0.0031 &  0.0917 & 0.0270 & 0.0008 &  0.1923 & 0.0567 & 0.0033 &  0.1884 & 0.0275 & 0.0009\\
      && & $\hat d_2$ &  0.3840 & 0.0618 & 0.0041 &  0.3855 & 0.0308 & 0.0012 &  0.2789 & 0.0594 & 0.0040 &  0.2849 & 0.0298 & 0.0011\\
      \cline{3-16}
      &&\multirow{2}{*}{ 0.9   }  &$\hat d_1$&0.0977 & 0.0573 & 0.0033 &  0.0931 & 0.0270 & 0.0008 &  0.1992 & 0.0578 & 0.0033 &  0.1891 & 0.0271 & 0.0009\\
      && & $\hat d_2$ &  0.3900 & 0.0607 & 0.0038 &  0.3815 & 0.0295 & 0.0012 &  0.2871 & 0.0599 & 0.0038 &  0.2844 & 0.0292 & 0.0011\\
      \cline{2-16}
      &\multirow{2}{*}{  \emph{LOB}   }& \multirow{2}{*}{  -   } &$\hat d_1$&0.1051 & 0.0579 & 0.0034 &  0.0955 & 0.0271 & 0.0008 &  0.2070 & 0.0580 & 0.0034 &  0.1915 & 0.0271 & 0.0008\\
&& & $\hat d_2$ &  0.3957 & 0.0606 & 0.0037 &  0.3822 & 0.0292 & 0.0012 &  0.2938 & 0.0602 & 0.0037 &  0.2860 & 0.0291 & 0.0010\\
\cline{2-16}
      &\multirow{2}{*}{  \emph{TLOB}   }& \multirow{2}{*}{  -   } &$\hat d_1$&0.1056 & 0.0785 & 0.0062 &  0.0956 & 0.0377 & 0.0014 &  0.2075 & 0.0783 & 0.0062 &  0.1916 & 0.0376 & 0.0015\\
    && & $\hat d_2$ &  0.4082 & 0.0773 & 0.0060 &  0.3866 & 0.0384 & 0.0017 &  0.2999 & 0.0773 & 0.0060 &  0.2882 & 0.0385 & 0.0016\\
      \hline\hline
      \multirow{12}{*}{$0.3$}&\multirow{4}{*}{  \emph{SLOB}   }&\multirow{2}{*}{ 0.7 }  &$\hat d_1$&0.1032 & 0.0534 & 0.0029 &  0.0974 & 0.0262 & 0.0007 &  0.2075 & 0.0556 & 0.0031 &  0.1947 & 0.0269 & 0.0008\\
&& & $\hat d_2$ &  0.4505 & 0.0776 & 0.0086 &  0.4162 & 0.0397 & 0.0018 &  0.3182 & 0.0630 & 0.0043 &  0.3003 & 0.0313 & 0.0010\\
\cline{3-16}
      &&\multirow{2}{*}{ 0.9   }  &$\hat d_1$&0.1044 & 0.0552 & 0.0031 &  0.0968 & 0.0261 & 0.0007 &  0.2075 & 0.0552 & 0.0031 &  0.1924 & 0.0262 & 0.0007\\
&& & $\hat d_2$ &  0.4261 & 0.0614 & 0.0045 &  0.3969 & 0.0300 & 0.0009 &  0.3095 & 0.0571 & 0.0033 &  0.2923 & 0.0281 & 0.0008\\
\cline{2-16}
      &\multirow{4}{*}{  \emph{SLOB}\hspace{-.25cm}$\phantom{B}^\ast$   }&\multirow{2}{*}{ 0.7 }  &$\hat d_1$&0.0910 & 0.0540 & 0.0030 &  0.0927 & 0.0263 & 0.0007 &  0.1904 & 0.0554 & 0.0032 &  0.1882 & 0.0267 & 0.0009\\
      && & $\hat d_2$ &  0.3878 & 0.0622 & 0.0040 &  0.3880 & 0.0304 & 0.0011 &  0.2811 & 0.0591 & 0.0038 &  0.2857 & 0.0291 & 0.0010\\
      \cline{3-16}
      &&\multirow{2}{*}{ 0.9   }  &$\hat d_1$&0.0965 & 0.0565 & 0.0032 &  0.0941 & 0.0264 & 0.0007 &  0.1975 & 0.0566 & 0.0032 &  0.1889 & 0.0264 & 0.0008\\
      && & $\hat d_2$ &  0.3928 & 0.0608 & 0.0037 &  0.3835 & 0.0288 & 0.0011 &  0.2890 & 0.0595 & 0.0037 &  0.2849 & 0.0284 & 0.0010\\
      \cline{2-16}
      &\multirow{2}{*}{  \emph{LOB}   }& \multirow{2}{*}{  -   } &$\hat d_1$&0.1036 & 0.0571 & 0.0033 &  0.0964 & 0.0264 & 0.0007 &  0.2052 & 0.0568 & 0.0032 &  0.1912 & 0.0264 & 0.0008\\
&& & $\hat d_2$ &  0.3983 & 0.0606 & 0.0037 &  0.3840 & 0.0285 & 0.0011 &  0.2958 & 0.0597 & 0.0036 &  0.2865 & 0.0282 & 0.0010\\
\cline{2-16}
      &\multirow{2}{*}{  \emph{TLOB}   }& \multirow{2}{*}{  -   } &$\hat d_1$&0.1062 & 0.0778 & 0.0061 &  0.0973 & 0.0372 & 0.0014 &  0.2071 & 0.0772 & 0.0060 &  0.1918 & 0.0369 & 0.0014\\
    && & $\hat d_2$ &  0.4106 & 0.0766 & 0.0060 &  0.3881 & 0.0374 & 0.0015 &  0.3010 & 0.0763 & 0.0058 &  0.2883 & 0.0374 & 0.0015\\
\hline\hline
      \multirow{12}{*}{$0.6$}&\multirow{4}{*}{  \emph{SLOB}   }&\multirow{2}{*}{ 0.7 }  &$\hat d_1$&0.1083 & 0.0509 & 0.0027 &  0.1064 & 0.0250 & 0.0007 &  0.2072 & 0.0517 & 0.0027 &  0.1959 & 0.0250 & 0.0006\\
&& & $\hat d_2$ &  0.4582 & 0.0752 & 0.0090 &  0.4243 & 0.0387 & 0.0021 &  0.3208 & 0.0589 & 0.0039 &  0.3014 & 0.0291 & 0.0008\\
\cline{3-16}
      &&\multirow{2}{*}{ 0.9   }  &$\hat d_1$&0.1071 & 0.0525 & 0.0028 &  0.1038 & 0.0247 & 0.0006 &  0.2064 & 0.0514 & 0.0027 &  0.1933 & 0.0244 & 0.0006\\
&& & $\hat d_2$ &  0.4322 & 0.0591 & 0.0045 &  0.4036 & 0.0285 & 0.0008 &  0.3120 & 0.0532 & 0.0030 &  0.2933 & 0.0259 & 0.0007\\
\cline{2-16}
      &\multirow{4}{*}{  \emph{SLOB}\hspace{-.25cm}$\phantom{B}^\ast$   }&\multirow{2}{*}{ 0.7 }  &$\hat d_1$&0.0918 & 0.0512 & 0.0027 &  0.0991 & 0.0248 & 0.0006 &  0.1877 & 0.0514 & 0.0028 &  0.1884 & 0.0248 & 0.0007\\
      && & $\hat d_2$ &  0.3967 & 0.0606 & 0.0037 &  0.3964 & 0.0291 & 0.0009 &  0.2852 & 0.0556 & 0.0033 &  0.2873 & 0.0269 & 0.0009\\
      \cline{3-16}
      &&\multirow{2}{*}{ 0.9   }  &$\hat d_1$&0.0966 & 0.0538 & 0.0029 &  0.0999 & 0.0250 & 0.0006 &  0.1950 & 0.0528 & 0.0028 &  0.1892 & 0.0246 & 0.0007\\
      && & $\hat d_2$ &  0.3995 & 0.0585 & 0.0034 &  0.3904 & 0.0270 & 0.0008 &  0.2924 & 0.0557 & 0.0032 &  0.2862 & 0.0261 & 0.0009\\
      \cline{2-16}
      &\multirow{2}{*}{  \emph{LOB}   }& \multirow{2}{*}{  -   } &$\hat d_1$&0.1032 & 0.0545 & 0.0030 &  0.1019 & 0.0250 & 0.0006 &  0.2025 & 0.0529 & 0.0028 &  0.1915 & 0.0246 & 0.0007\\
&& & $\hat d_2$ &  0.4045 & 0.0580 & 0.0034 &  0.3906 & 0.0266 & 0.0008 &  0.2991 & 0.0557 & 0.0031 &  0.2878 & 0.0260 & 0.0008\\
\cline{2-16}
      &\multirow{2}{*}{  \emph{TLOB}   }& \multirow{2}{*}{  -   } &$\hat d_1$&0.1087 & 0.0742 & 0.0056 &  0.1040 & 0.0356 & 0.0013 &  0.2060 & 0.0718 & 0.0052 &  0.1927 & 0.0346 & 0.0012\\
    && & $\hat d_2$ &  0.4164 & 0.0735 & 0.0057 &  0.3943 & 0.0352 & 0.0013 &  0.3035 & 0.0713 & 0.0051 &  0.2890 & 0.0345 & 0.0013\\
\hline\hline
      \multirow{12}{*}{$0.8$}&\multirow{4}{*}{  \emph{SLOB}   }&\multirow{2}{*}{ 0.7 }  &$\hat d_1$&0.1189 & 0.0494 & 0.0028 &  0.1252 & 0.0253 & 0.0013 &  0.2089 & 0.0476 & 0.0023 &  0.1990 & 0.0234 & 0.0005\\
&& & $\hat d_2$ &  0.4697 & 0.0740 & 0.0103 &  0.4411 & 0.0388 & 0.0032 &  0.3241 & 0.0551 & 0.0036 &  0.3041 & 0.0271 & 0.0008\\
\cline{3-16}
      &&\multirow{2}{*}{ 0.9   }  &$\hat d_1$&0.1142 & 0.0506 & 0.0028 &  0.1193 & 0.0241 & 0.0010 &  0.2068 & 0.0473 & 0.0023 &  0.1956 & 0.0226 & 0.0005\\
&& & $\hat d_2$ &  0.4409 & 0.0581 & 0.0050 &  0.4178 & 0.0278 & 0.0011 &  0.3147 & 0.0494 & 0.0026 &  0.2954 & 0.0238 & 0.0006\\
\cline{2-16}
      &\multirow{4}{*}{  \emph{SLOB}\hspace{-.25cm}$\phantom{B}^\ast$   }&\multirow{2}{*}{ 0.7 }  &$\hat d_1$&0.0976 & 0.0495 & 0.0024 &  0.1146 & 0.0240 & 0.0008 &  0.1868 & 0.0473 & 0.0024 &  0.1904 & 0.0230 & 0.0006\\
      && & $\hat d_2$ &  0.4088 & 0.0594 & 0.0036 &  0.4128 & 0.0284 & 0.0010 &  0.2893 & 0.0517 & 0.0028 &  0.2900 & 0.0248 & 0.0007\\
      \cline{3-16}
      &&\multirow{2}{*}{ 0.9   }  &$\hat d_1$&0.1007 & 0.0520 & 0.0027 &  0.1138 & 0.0240 & 0.0008 &  0.1939 & 0.0487 & 0.0024 &  0.1910 & 0.0229 & 0.0006\\
      && & $\hat d_2$ &  0.4086 & 0.0566 & 0.0033 &  0.4046 & 0.0257 & 0.0007 &  0.2956 & 0.0514 & 0.0027 &  0.2885 & 0.0239 & 0.0007\\
      \cline{2-16}
      &\multirow{2}{*}{  \emph{LOB}   }& \multirow{2}{*}{  -   } &$\hat d_1$&0.1064 & 0.0528 & 0.0028 &  0.1152 & 0.0241 & 0.0008 &  0.2011 & 0.0490 & 0.0024 &  0.1931 & 0.0228 & 0.0006\\
&& & $\hat d_2$ &  0.4125 & 0.0559 & 0.0033 &  0.4042 & 0.0252 & 0.0007 &  0.3022 & 0.0512 & 0.0026 &  0.2899 & 0.0237 & 0.0007\\
\cline{2-16}
      &\multirow{2}{*}{  \emph{TLOB}   }& \multirow{2}{*}{  -   } &$\hat d_1$&0.1145 & 0.0715 & 0.0053 &  0.1184 & 0.0344 & 0.0015 &  0.2058 & 0.0656 & 0.0043 &  0.1948 & 0.0321 & 0.0011\\
    && & $\hat d_2$ &  0.4247 & 0.0708 & 0.0056 &  0.4078 & 0.0337 & 0.0012 &  0.3067 & 0.0653 & 0.0043 &  0.2911 & 0.0316 & 0.0011\\
\hline \hline
\end{tabular}}
\end{table}

\begin{table}[!h]
\renewcommand{\arraystretch}{1.2}
\setlength{\tabcolsep}{3pt}
\caption{Simulation results for the estimators \emph{SLOB}, $SLOB^\ast$, \emph{TLOB} and \emph{LOB} in Gaussian VARFIMA$(0,\bs d,0)$ processes. Presented are the estimated values (mean), their standard deviation (st.d) and mean square error of the estimates (mse).}\label{tab+2}
\vskip.3cm
\centering
{\scriptsize
\begin{tabular}{|c|c|c|c|ccc|ccc|ccc|ccc|}
\hline \hline
\multirow{3}{*}{$\rho$}&\multirow{3}{*}{Method}&\multirow{3}{*}{$\beta$}&   \multirow{3}{*}{$\hat d_i$}    &\multicolumn{6}{|c|}{$\bs d=(0.1,0.3)$}
      &\multicolumn{6}{c|}{$\bs d=(0.3,0.4)$}   \\
      \cline{5-16}
      &       &       &      &\multicolumn{3}{|c|}{$\alpha=0.65$}       &\multicolumn{3}{|c|}{$\alpha=0.85$}            &\multicolumn{3}{|c|}{$\alpha=0.65$}
      &\multicolumn{3}{|c|}{$\alpha=0.85$}    \\
      \cline{5-16}
      &       &       &      & mean & st.d. & mse& mean & st.d. & mse& mean & st.d. & mse& mean & st.d. & mse\\
      \hline\hline
      \multirow{12}{*}{$0$}&\multirow{4}{*}{  \emph{SLOB}   }&\multirow{2}{*}{ 0.7 }  &$\hat d_1$&0.1024 & 0.0543 & 0.0030 &  0.0954 & 0.0268 & 0.0007 &  0.3228 & 0.0614 & 0.0043 &  0.2986 & 0.0299 & 0.0009\\
&& & $\hat d_2$ &  0.3168 & 0.0644 & 0.0044 &  0.2999 & 0.0321 & 0.0010 &  0.4462 & 0.0786 & 0.0083 &  0.4137 & 0.0401 & 0.0018\\
\cline{3-16}
      &&\multirow{2}{*}{ 0.9   }  &$\hat d_1$&0.1044 & 0.0560 & 0.0032 &  0.0953 & 0.0267 & 0.0007 &  0.3166 & 0.0573 & 0.0036 &  0.2915 & 0.0274 & 0.0008\\
&& & $\hat d_2$ &  0.3080 & 0.0580 & 0.0034 &  0.2919 & 0.0289 & 0.0009 &  0.4232 & 0.0624 & 0.0044 &  0.3949 & 0.0306 & 0.0010\\
\cline{2-16}
      &\multirow{4}{*}{  \emph{SLOB}\hspace{-.25cm}$\phantom{B}^\ast$   }&\multirow{2}{*}{ 0.7 }  &$\hat d_1$&0.0922 & 0.0549 & 0.0031 &  0.0917 & 0.0270 & 0.0008 &  0.2969 & 0.0590 & 0.0035 &  0.2879 & 0.0285 & 0.0010\\
      && & $\hat d_2$ &  0.2790 & 0.0594 & 0.0040 &  0.2849 & 0.0298 & 0.0011 &  0.3838 & 0.0617 & 0.0041 &  0.3855 & 0.0308 & 0.0012\\
      \cline{3-16}
      &&\multirow{2}{*}{ 0.9   }  &$\hat d_1$&0.0977 & 0.0573 & 0.0033 &  0.0931 & 0.0270 & 0.0008 &  0.3026 & 0.0585 & 0.0034 &  0.2862 & 0.0274 & 0.0009\\
      && & $\hat d_2$ &  0.2872 & 0.0599 & 0.0038 &  0.2844 & 0.0292 & 0.0011 &  0.3899 & 0.0607 & 0.0038 &  0.3815 & 0.0295 & 0.0012\\
      \cline{2-16}
      &\multirow{2}{*}{  \emph{LOB}   }& \multirow{2}{*}{  -   } &$\hat d_1$&0.1051 & 0.0580 & 0.0034 &  0.0955 & 0.0271 & 0.0008 &  0.3098 & 0.0582 & 0.0035 &  0.2880 & 0.0272 & 0.0009\\
&& & $\hat d_2$ &  0.2938 & 0.0602 & 0.0037 &  0.2860 & 0.0291 & 0.0010 &  0.3957 & 0.0605 & 0.0037 &  0.3822 & 0.0292 & 0.0012\\
\cline{2-16}
      &\multirow{2}{*}{  \emph{TLOB}   }& \multirow{2}{*}{  -   } &$\hat d_1$&0.1055 & 0.0785 & 0.0062 &  0.0956 & 0.0377 & 0.0014 &  0.3112 & 0.0780 & 0.0062 &  0.2882 & 0.0375 & 0.0015\\
    && & $\hat d_2$ &  0.2999 & 0.0773 & 0.0060 &  0.2882 & 0.0385 & 0.0016 &  0.4083 & 0.0774 & 0.0060 &  0.3866 & 0.0385 & 0.0017\\
\hline \hline
      \multirow{12}{*}{$0.3$}&\multirow{4}{*}{  \emph{SLOB}   }&\multirow{2}{*}{ 0.7 }  &$\hat d_1$&0.1024 & 0.0533 & 0.0028 &  0.0418 & 0.0248 & 0.0040 &  0.3219 & 0.0599 & 0.0041 &  0.2988 & 0.0290 & 0.0008\\
&& & $\hat d_2$ &  0.3190 & 0.0633 & 0.0044 &  0.2345 & 0.0319 & 0.0053 &  0.4482 & 0.0766 & 0.0082 &  0.4144 & 0.0393 & 0.0017\\
\cline{3-16}
      &&\multirow{2}{*}{ 0.9   }  &$\hat d_1$&0.1038 & 0.0550 & 0.0030 &  0.0959 & 0.0260 & 0.0007 &  0.3155 & 0.0560 & 0.0034 &  0.2915 & 0.0266 & 0.0008\\
&& & $\hat d_2$ &  0.3099 & 0.0572 & 0.0034 &  0.2928 & 0.0281 & 0.0008 &  0.4250 & 0.0610 & 0.0043 &  0.3954 & 0.0298 & 0.0009\\
\cline{2-16}
      &\multirow{4}{*}{  \emph{SLOB}\hspace{-.25cm}$\phantom{B}^\ast$   }&\multirow{2}{*}{ 0.7 }  &$\hat d_1$&0.0907 & 0.0539 & 0.0030 &  0.0368 & 0.0249 & 0.0046 &  0.2946 & 0.0575 & 0.0033 &  0.2874 & 0.0276 & 0.0009\\
      && & $\hat d_2$ &  0.2817 & 0.0594 & 0.0039 &  0.2163 & 0.0278 & 0.0078 &  0.3864 & 0.0614 & 0.0039 &  0.3864 & 0.0302 & 0.0011\\
      \cline{3-16}
      &&\multirow{2}{*}{ 0.9   }  &$\hat d_1$&0.0962 & 0.0564 & 0.0032 &  0.0934 & 0.0263 & 0.0007 &  0.3007 & 0.0572 & 0.0033 &  0.2859 & 0.0267 & 0.0009\\
      && & $\hat d_2$ &  0.2893 & 0.0597 & 0.0037 &  0.2854 & 0.0284 & 0.0010 &  0.3920 & 0.0602 & 0.0037 &  0.3822 & 0.0287 & 0.0011\\
      \cline{2-16}
      &\multirow{2}{*}{  \emph{LOB}   }& \multirow{2}{*}{  -   } &$\hat d_1$&0.1034 & 0.0570 & 0.0033 &  0.0957 & 0.0264 & 0.0007 &  0.3078 & 0.0570 & 0.0033 &  0.2877 & 0.0265 & 0.0009\\
&& & $\hat d_2$ &  0.2960 & 0.0599 & 0.0036 &  0.2870 & 0.0283 & 0.0010 &  0.3978 & 0.0600 & 0.0036 &  0.3828 & 0.0284 & 0.0011\\
\cline{2-16}
      &\multirow{2}{*}{  \emph{TLOB}   }& \multirow{2}{*}{  -   } &$\hat d_1$&0.1055 & 0.0775 & 0.0060 &  0.0964 & 0.0371 & 0.0014 &  0.3108 & 0.0770 & 0.0060 &  0.2885 & 0.0368 & 0.0015\\
    && & $\hat d_2$ &  0.3014 & 0.0765 & 0.0059 &  0.2888 & 0.0374 & 0.0015 &  0.4096 & 0.0761 & 0.0059 &  0.3868 & 0.0373 & 0.0016\\
\hline \hline
      \multirow{12}{*}{$0.6$}&\multirow{4}{*}{  \emph{SLOB}   }&\multirow{2}{*}{ 0.7 }  &$\hat d_1$&0.1043 & 0.0500 & 0.0025 &  0.0045 & 0.0231 & 0.0096 &  0.3221 & 0.0558 & 0.0036 &  0.3004 & 0.0270 & 0.0007\\
&& & $\hat d_2$ &  0.3237 & 0.0602 & 0.0042 &  0.1921 & 0.0319 & 0.0127 &  0.4510 & 0.0718 & 0.0077 &  0.4158 & 0.0371 & 0.0016\\
\cline{3-16}
      &&\multirow{2}{*}{ 0.9   }  &$\hat d_1$&0.1042 & 0.0517 & 0.0027 &  0.0991 & 0.0244 & 0.0006 &  0.3147 & 0.0521 & 0.0029 &  0.2925 & 0.0247 & 0.0007\\
&& & $\hat d_2$ &  0.3138 & 0.0541 & 0.0031 &  0.2961 & 0.0261 & 0.0007 &  0.4277 & 0.0570 & 0.0040 &  0.3966 & 0.0277 & 0.0008\\
\cline{2-16}
      &\multirow{4}{*}{  \emph{SLOB}\hspace{-.25cm}$\phantom{B}^\ast$   }&\multirow{2}{*}{ 0.7 }  &$\hat d_1$&0.0898 & 0.0505 & 0.0027 &  -0.0025 & 0.0231 & 0.0110 &  0.2913 & 0.0530 & 0.0029 &  0.2874 & 0.0254 & 0.0008\\
      && & $\hat d_2$ &  0.2874 & 0.0567 & 0.0034 &  0.1719 & 0.0268 & 0.0171 &  0.3912 & 0.0579 & 0.0034 &  0.3885 & 0.0282 & 0.0009\\
      \cline{3-16}
      &&\multirow{2}{*}{ 0.9   }  &$\hat d_1$&0.0949 & 0.0531 & 0.0028 &  0.0960 & 0.0247 & 0.0006 &  0.2978 & 0.0532 & 0.0028 &  0.2860 & 0.0248 & 0.0008\\
      && & $\hat d_2$ &  0.2938 & 0.0565 & 0.0032 &  0.2889 & 0.0263 & 0.0008 &  0.3958 & 0.0564 & 0.0032 &  0.3836 & 0.0265 & 0.0010\\
      \cline{2-16}
      &\multirow{2}{*}{  \emph{LOB}   }& \multirow{2}{*}{  -   } &$\hat d_1$&0.1017 & 0.0537 & 0.0029 &  0.0981 & 0.0247 & 0.0006 &  0.3049 & 0.0530 & 0.0028 &  0.2878 & 0.0246 & 0.0008\\
&& & $\hat d_2$ &  0.3003 & 0.0565 & 0.0032 &  0.2903 & 0.0261 & 0.0008 &  0.4014 & 0.0561 & 0.0031 &  0.3842 & 0.0261 & 0.0009\\
\cline{2-16}
      &\multirow{2}{*}{  \emph{TLOB}   }& \multirow{2}{*}{  -   } &$\hat d_1$&0.1058 & 0.0729 & 0.0053 &  0.0996 & 0.0351 & 0.0012 &  0.3100 & 0.0718 & 0.0052 &  0.2895 & 0.0345 & 0.0013\\
    && & $\hat d_2$ &  0.3053 & 0.0722 & 0.0052 &  0.2917 & 0.0348 & 0.0013 &  0.4122 & 0.0713 & 0.0052 &  0.3876 & 0.0345 & 0.0013\\
\hline \hline
      \multirow{12}{*}{$0.8$}&\multirow{4}{*}{  \emph{SLOB}   }&\multirow{2}{*}{ 0.7 }  &$\hat d_1$&0.1103 & 0.0473 & 0.0023 &  -0.0111 & 0.0222 & 0.0128 &  0.3247 & 0.0517 & 0.0033 &  0.3045 & 0.0254 & 0.0007\\
&& & $\hat d_2$ &  0.3312 & 0.0583 & 0.0044 &  0.1752 & 0.0323 & 0.0166 &  0.4546 & 0.0671 & 0.0075 &  0.4190 & 0.0352 & 0.0016\\
\cline{3-16}
      &&\multirow{2}{*}{ 0.9   }  &$\hat d_1$&0.1076 & 0.0487 & 0.0024 &  0.1070 & 0.0232 & 0.0006 &  0.3157 & 0.0480 & 0.0026 &  0.2952 & 0.0230 & 0.0006\\
&& & $\hat d_2$ &  0.3194 & 0.0516 & 0.0030 &  0.3035 & 0.0245 & 0.0006 &  0.4306 & 0.0530 & 0.0037 &  0.3990 & 0.0258 & 0.0007\\
\cline{2-16}
      &\multirow{4}{*}{  \emph{SLOB}\hspace{-.25cm}$\phantom{B}^\ast$   }&\multirow{2}{*}{ 0.7 }  &$\hat d_1$&0.0925 & 0.0477 & 0.0023 &  -0.0202 & 0.0220 & 0.0149 &  0.2898 & 0.0487 & 0.0025 &  0.2892 & 0.0235 & 0.0007\\
      && & $\hat d_2$ &  0.2953 & 0.0542 & 0.0030 &  0.1536 & 0.0263 & 0.0221 &  0.3958 & 0.0540 & 0.0029 &  0.3918 & 0.0262 & 0.0008\\
      \cline{3-16}
      &&\multirow{2}{*}{ 0.9   }  &$\hat d_1$&0.0965 & 0.0502 & 0.0025 &  0.1031 & 0.0234 & 0.0006 &  0.2964 & 0.0490 & 0.0024 &  0.2877 & 0.0230 & 0.0007\\
      && & $\hat d_2$ &  0.2997 & 0.0534 & 0.0028 &  0.2964 & 0.0244 & 0.0006 &  0.3993 & 0.0521 & 0.0027 &  0.3861 & 0.0243 & 0.0008\\
      \cline{2-16}
      &\multirow{2}{*}{  \emph{LOB}   }& \multirow{2}{*}{  -   } &$\hat d_1$&0.1026 & 0.0509 & 0.0026 &  0.1049 & 0.0234 & 0.0006 &  0.3033 & 0.0490 & 0.0024 &  0.2893 & 0.0228 & 0.0006\\
&& & $\hat d_2$ &  0.3056 & 0.0532 & 0.0029 &  0.2975 & 0.0242 & 0.0006 &  0.4047 & 0.0517 & 0.0027 &  0.3865 & 0.0239 & 0.0007\\
\cline{2-16}
      &\multirow{2}{*}{  \emph{TLOB}   }& \multirow{2}{*}{  -   } &$\hat d_1$&0.1083 & 0.0684 & 0.0047 &  0.1069 & 0.0332 & 0.0011 &  0.3104 & 0.0657 & 0.0044 &  0.2919 & 0.0320 & 0.0011\\
    && & $\hat d_2$ &  0.3110 & 0.0678 & 0.0047 &  0.2989 & 0.0325 & 0.0011 &  0.4155 & 0.0655 & 0.0045 &  0.3898 & 0.0317 & 0.0011\\
\hline\hline
\end{tabular}}
\end{table}

\FloatBarrier
Overall, all estimators perform well, but the \emph{SLOB} estimator usually presents better performance in terms of mse and st.d. Usually the estimator with the smallest mse is not the one with the smallest bias. The \emph{TLOB} is the estimator with the worst performance in terms of mse. Estimator \emph{SLOB} outperforms the estimator $SLOB^\ast$ in all cases, except for $\bs d=(0.1,0.4)$ when $\rho\in\{0.6,0.8\}$ where the latter yields sensibly better estimates than the former. Overall, the best estimator in terms of mse is the \emph{SLOB} with $(\alpha,\beta)=(0.85,0.9)$.

As for bias, for the \emph{LOB} and \emph{TLOB} estimators, $\alpha = 0.65$ yields estimates with smallest bias\footnote{measured as the sum of the absolute distance between the estimates $(\widehat{d}_1,\widehat{d}_2)$ and $\do$.} in most of cases (14 and 10 out of 16 cases, respectively). For the $SLOB^\ast$ estimator, the combination $(\alpha,\beta)=(0.65,0.9)$ is the one presenting smallest bias in most cases (14 out of 16 cases), while for the \emph{SLOB} estimator $(\alpha,\beta)=(0.85,0.9)$ and $(0.65,0.9)$ are the combinations yielding the smallest bias, with a small advantage to the former (9 and 7, out of 16 cases, respectively).

\begin{figure}[h!]
\centering
\mbox{
 \subfigure[]{\includegraphics[width=0.2\textwidth]{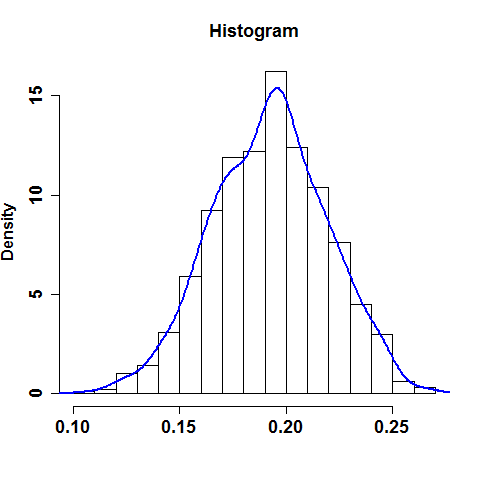}}\hskip.2cm
 \subfigure[]{\includegraphics[width=0.2\textwidth]{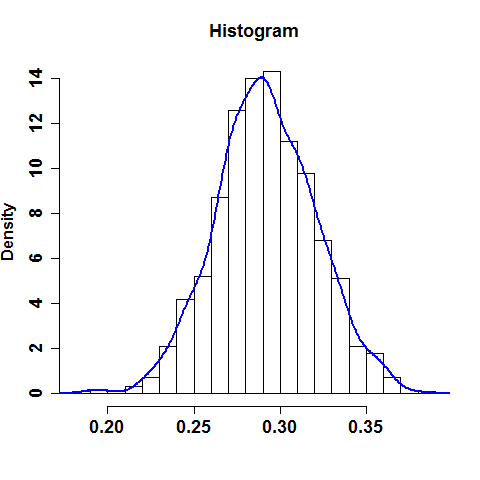}}\hskip.2cm
 \subfigure[]{\includegraphics[width=0.2\textwidth]{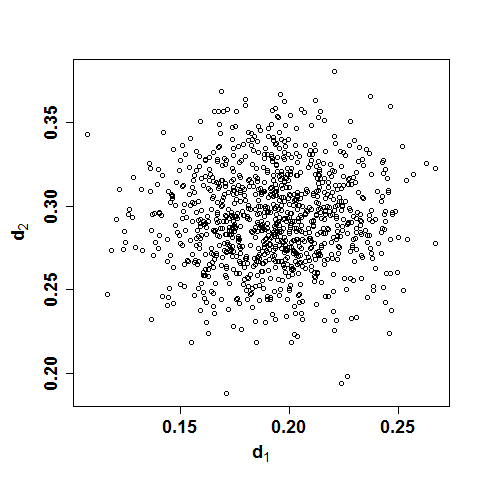}}
 }
 \mbox{
 \subfigure[]{\includegraphics[width=0.2\textwidth]{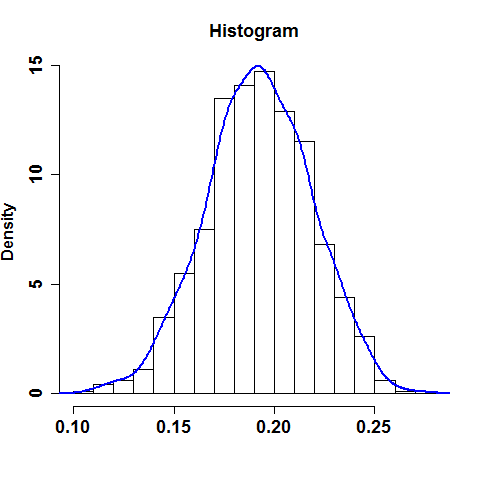}}\hskip.2cm
 \subfigure[]{\includegraphics[width=0.2\textwidth]{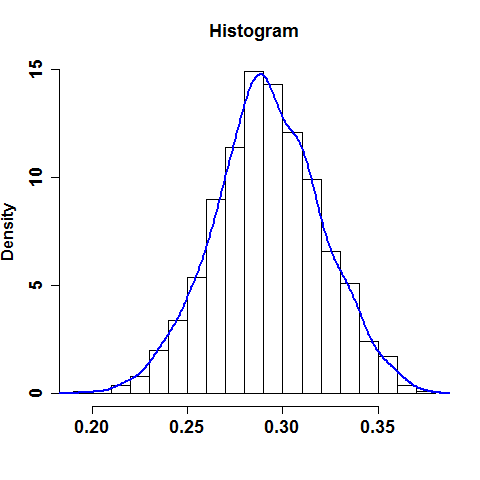}}\hskip.2cm
 \subfigure[]{\includegraphics[width=0.2\textwidth]{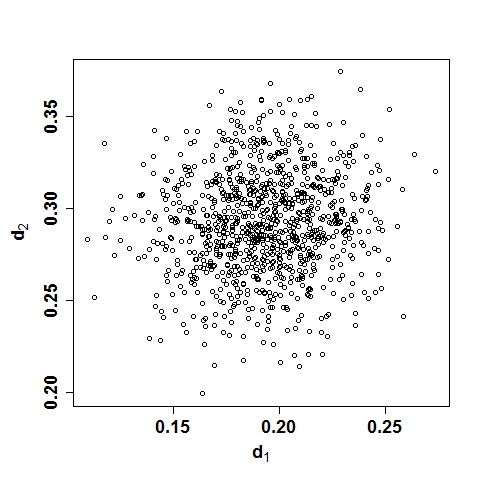}}
  }
  \mbox{
 \subfigure[]{\includegraphics[width=0.2\textwidth]{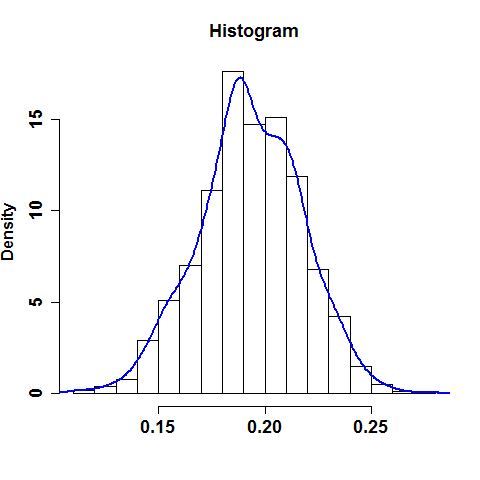}}\hskip.2cm
 \subfigure[]{\includegraphics[width=0.2\textwidth]{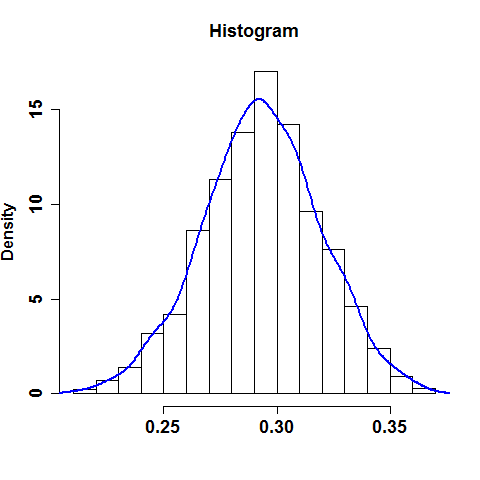}}\hskip.2cm
  \subfigure[]{\includegraphics[width=0.2\textwidth]{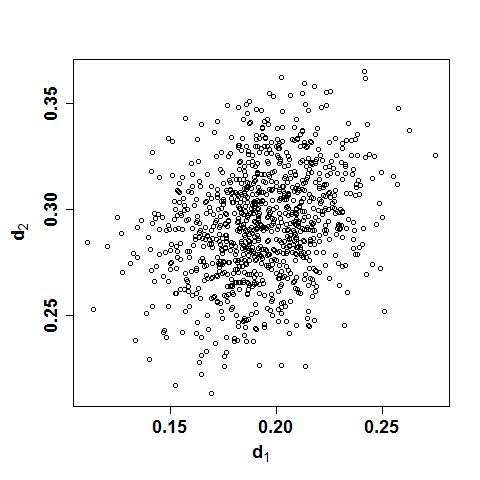}}
  }
  \mbox{
 \subfigure[]{\includegraphics[width=0.2\textwidth]{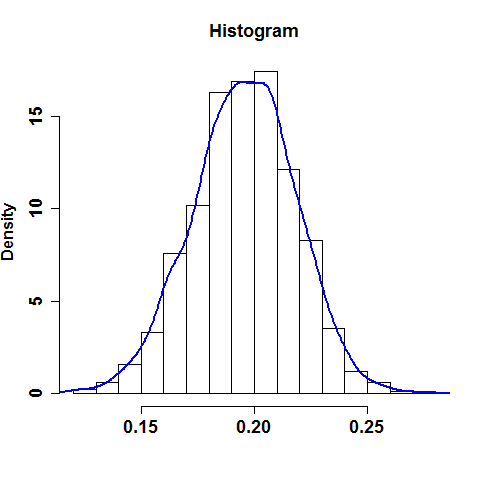}}\hskip.2cm
 \subfigure[]{\includegraphics[width=0.2\textwidth]{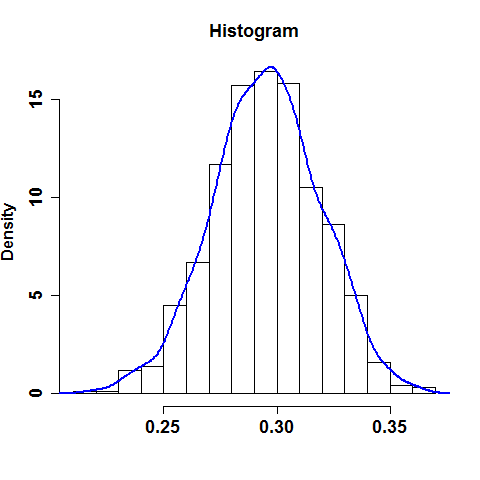}}\hskip.2cm
   \subfigure[]{\includegraphics[width=0.2\textwidth]{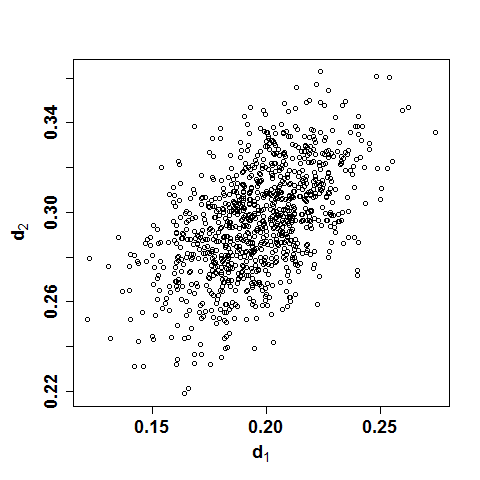}}
  }
  \caption{Histogram, kernel density and scatter plot of the \emph{SLOB} estimated values of $\bs d_0=(0.2,0.3)$ for (a)--(c) $\rho =0$; (d)--(f) $\rho =0.3$; (g)--(i) $\rho =0.6$ and (j)--(l) $\rho=0.8$. }\label{hist}
\end{figure}
\FloatBarrier

We also notice that the trade-off between bias and variance present in estimators of the spectral density function based on the smoothed periodogram does not seem to influence the estimation of the parameter $\bs d$. This is expected since the smoothed periodogram is a function of the observed time series alone which, by its turn, depends only on $\do$ and not on any particular estimated value $\widehat{\bs d}$, being, hence, a constant with respect to the latter.

Figure \ref{hist}  depicts the behavior of the estimated values by presenting the histogram, kernel density estimates and the scatter plot of the \emph{SLOB} estimated values for $(\alpha,\beta)=(0.85,0.9)$ and $\bs d =(0.2,0.3)$. Figures \ref{hist}(a)--(c) correspond to $\rho=0$, Figures \ref{hist}(d)--(f) to $\rho=0.3$, Figures \ref{hist}(g)--(i) to $\rho=0.6$ and Figures \ref{hist}(j)--(l) to $\rho=0.8$. Figure \ref{hist} reinforces the conjecture that the \emph{SLOB} estimator is indeed asymptotically normally distributed.

\section{Empirical Application}

In this section we apply the \emph{LOB}, \emph{SLOB} and \emph{TLOB} estimators considered in the last section to an exchange rate data set consisting of four daily exchange rates (business days) against the British Pound, namely, the US Dollar (USD/GBP), the Euro (EUR/GBP), the Japanese Yen (JPY/GBP) and the Swiss Franc (CHF/GBP). The data comprehend the period between October 2008 and July 2012, with sample size $n=1,684$. The data is similar to the one applied in Lobato (1999). The time series considered are the squared log-returns (squares of the first difference of the logarithm of the exchange rates) which is usually associated to the volatility of the returns (see, for instance, Cont, 2001). Our goal is to estimate the fractional differencing parameter $\bs d=(d_1,d_2,d_3,d_4)^\prime$, where the components are associated to the USD/GBP, the EUR/GBP, the JPY/GBP and the CHF/GBP exchange rates, respectively.   In accordance to our findings in Section \ref{sim}, for all estimators we apply $m=\lfloor n^{0.85}\rfloor=552$, while the cut-off point of the smoothed periodogram needed for the \emph{SLOB} estimator is $\ell(n)=\lfloor n^{0.9}\rfloor=801$. For the \emph{SLOB} estimator we apply the Bartlett's weights and for the \emph{TLOB} estimator, we apply the cosine-bell tapering function. For the \emph{SLOB} estimator, we are assuming that it is indeed asymptotically normally distributed although, at this moment, we could not present a formal proof of the result, but we have empirical evidences that this is indeed true (see Table \ref{tab+1}, Table \ref{tab+2} and Figure \ref{hist}).

The correlation matrix of the data is presented on Table \ref{etab}. As expected, all correlations are positive and high. Worth of note is the very high correlation between the squared returns of the USD/GBP and JPY/GBP exchange rates, over 0.85, which suggests a high association between the volatility on their exchange rates against the British Pound. Table \ref{etab} also presents the estimates of the fractional differencing parameter $\bs d$ obtained. In all cases the estimators pointed toward the existence of a mild long-range dependence on the data, characterized by the relatively small positive values of the estimated fractional differencing parameter. The estimated values obtained from the \emph{SLOB} and \emph{TLOB} estimators were all higher than the ones obtained from the \emph{LOB} estimator.

\begin{table}[h]
\setlength{\tabcolsep}{3pt}
\renewcommand{\arraystretch}{1.2}
\caption{Estimated values of the fractional differencing parameter $\bs d$ and correlation matrix of the squared log-returns of the daily exchange rates data set for the \emph{LOB}, \emph{SLOB} and \emph{TLOB} estimators.}\label{etab}
\vskip.3cm
\centering
\footnotesize
\begin{tabular}{|c|c|c|c|cccc|}
\hline\hline
\multirow{2}{*}{Exchange rate} & \multirow{2}{*}{\emph{LOB}} &\multirow{2}{*}{\emph{SLOB}} & \multirow{2}{*}{\emph{TLOB}} & \multicolumn{4}{c|}{Correlation Matrix}\\
\cline{5-8}
&&&& USD/GBP&EUR/GBP&JPY/GBP&CHF/GBP\\
\hline\hline
USD/GBP & 0.1065 & 0.1584 & 0.1546 &      1      & 0.4534 & 0.8552 & 0.4404 \\
EUR/GBP & 0.1459 & 0.1960 & 0.1895 &      0.4534 &  1     & 0.3273 & 0.4495 \\
JPY/GBP & 0.0596 & 0.1040 & 0.1030 &      0.8552 & 0.3273 & 1      & 0.4272 \\
CHF/GBP & 0.0960 & 0.1450 & 0.1421 &      0.4404 & 0.4495 & 0.4272 & 1      \\
\hline\hline
\end{tabular}
\end{table} 
\normalsize

In both cases the estimated values point out for a greater persistence on the volatility of the EUR/GBP, USD/GBP, CHF/GBP and JPY/BPN exchange rates in this order. We observe that the estimated values of the fractional differencing parameter for the squared log-returns of the USD/GBP and the CHF/GBP exchange rates are remarkably close to one another. In order to apply the asymptotic results of Section \ref{asym}, we shall test for the equality between them as well as the equality among the four fractional differencing parameters.

In order to do that, let {\small $\widehat {\mathrm G}_L$, $\widehat {\mathrm G}_S$} and {\small$\widehat {\mathrm G}_T$} denote the estimators of $G_0$ obtained by taking $f_n$ as the periodogram of the series, the smoothed periodogram with the Bartlett's weights and the tapered periodogram with the cosine-bell tapering function, respectively. Also denote by {\small$\widehat\Omega_L$, $\widehat\Omega_S$} and {\small$\widehat\Omega_T$} the respective estimates of the matrix $\Omega$ as defined in Theorem \ref{norm} based on the pairs of estimates of $\d$ and $\widehat G$ obtained with the \emph{LOB}, \emph{SLOB} and \emph{TLOB} estimators and the respective estimates of $G_0$. The estimates obtained are

\footnotesize\[
\widehat\Omega_L=
\left(
 \begin{array}{cccc}
0.1410 & 0.0295 & 0.0776 & 0.0051\\
0.0295 & 0.2002 & 0.0074 & 0.0133\\
0.0776 & 0.0074 & 0.1502 & 0.0177\\
0.0051 & 0.0133 & 0.0177 & 0.2146\\
 \end{array}
\right),\qquad
\widehat\Omega_S=
\left(
 \begin{array}{cccc}
0.1438 & 0.0134 & 0.0867 & 0.0109\\
0.0134 & 0.2131 & 0.0042 & 0.0198\\
0.0867 & 0.0042 & 0.1505 & 0.0134\\
0.0109 & 0.0198 & 0.0134 & 0.2068\\
\end{array}
\right)
\]\normalsize
and
\footnotesize\[
\widehat\Omega_T=
\left(
 \begin{array}{cccc}
0.2495&  0.0003&  0.0000&  0.0002\\ 
0.0003&  0.2496&  0.0001&  0.0000\\
0.0000&  0.0001&  0.2497&  0.0002\\
0.0002&  0.0000&  0.0002&  0.2495\\
\end{array}
\right).\]\normalsize

Let us start by testing the equality of the fractional differencing parameters for the squared log-returns of the USD/GBP and the CHF/GBP exchange rates. In order words, we want to test the null hypothesis $H_0:d_1=d_4$. In this case, $R=(1,0,0,-1)$ and $\bs\nu=0$. The test statistics \eqref{test} under the null hypothesis is distributed according to a $\chi^2_1$ distribution. In this case, for the test statistics based on the \emph{LOB}, \emph{SLOB} and \emph{TLOB} estimators, we have $T_{L}=0.1838$, $T_{S}=0.2986$ and $T_T=0.1745$ with all $p$-values $>0.58$. Therefore, we cannot reject the null hypothesis at any reasonable significance level and we conclude that the fractional differencing parameter of squared log-returns of the USD/GBP and CHF/GBP exchange rates are statistically equal. The question that naturally arises is are all fractional differencing parameters statistically equal? That is, we are interested in testing $H_0:d_1=\cdots =d_4$. In this case

\footnotesize\[R=\left(
\begin{array}{cccc}
  1 & -1 & 0 & 0 \\
  0 & 1 & -1 & 0 \\
  0 & 0 & 1 & -1 \\
\end{array}\right)\qquad\mbox{and}\qquad \bs\nu=(0,0,0).
\]\normalsize

The values of the test statistics in this situation are $T_L=16.38$, $T_S=19.35$ and $T_T=8.45$, which, under the null hypothesis, are distributed according to a $\chi^2_3$ distribution. In this case, for the $T_L$ and $T_S$ test statistics, the $p$-values are both smaller than $0.001$ strongly indicating that the fractional differencing parameters of the series cannot be considered all equals at any reasonable significance level. The $T_T$ statistics also agrees (at 5\% confidence level) with this conclusion, but with $p$-value equals to $0.038$.

\section{Conclusions}
In this work we propose and study a class of Gaussian semiparametric estimator of the fractional differencing parameter in multivariate long-range dependent processes. The main idea is to modify the approach of Shimotsu (2007) in order to introduce and analyze a generalization of the two-step estimator proposed by Lobato (1999) and studied in a slightly different form by Shimotsu (2007). More specifically, the idea is to consider the same objective function as in Lobato (1999), but considering an arbitrary spectral density matrix estimator in place of the periodogram. We consider two main cases: first, for an arbitrary consistent spectral density estimator, we show that the proposed estimator is also consistent. Secondly, we consider a class of spectral density estimator satisfying a single condition and show that the resulting estimator is consistent as well. The asymptotic distribution of the estimator under a certain condition on the spectral density estimator is derived and shown to be Gaussian with the same variance-covariance matrix as the one derived in Shimotsu (2007). Hypothesis testing are also discussed in connection with the asymptotic theory.

To assess the finite sample performance of the proposed estimator and to compare to the original one, we present a Monte Carlo simulation study based on VARFIMA$(0,\bs d,0)$ processes. Under the conditions of the experiment, the proposed estimator shows an overall better performance over the original one. The approach is also applied to a financial data set consisting of daily exchange rates against the British Pound.


\subsubsection*{Acknowledgements}
\footnotesize
G. Pumi's research was supported by a REUNI Pos-Doctoral fellowship from CAPES-Brazil via the REUNI program. S.R.C. Lopes' research was partially supported by CNPq-Brazil, by Pronex {\it Probabilidade e Processos Estoc\'asticos} - E-26/170.008/2008 -APQ1 and also by INCT {\it em Matem\'atica}. The authors are grateful to the (Brazilian) National Center of Super Computing (CESUP-UFRGS) for the computational resources.

\normalsize

\section*{Appendix A: Proofs}
\renewcommand{\theequation}{A.\arabic{equation}}
\renewcommand{\thesubsection}{A}
\renewcommand{\theequation}{A.\arabic{equation}}
\setcounter{equation}{0}
In this section we present the proofs of the results in Sections 4 and 5. We establish lemmas and  theorems in the same sequence as they appear in the text.

\subsubsection{Proof of Theorem \ref{consistency}:}

The proof follows the same lines as the proof of theorem 3(a) in Shimotsu (2007). First, without loss of generality let $0<\delta<1/2$ be fixed and consider the set $N_\delta\vcentcolon=\big\{\bs d:\|\bs d-\do\|_{\infty}>\delta\big\}$. Define $\bs\te=(\te_1,\cdots,\te_q)^\prime\vcentcolon=\bs d-\do$ and let \small\[L(\bs d)\vcentcolon=S(\bs d)-S(\do).\]\normalsize
Let $0<\epsilon<1/4$ and define $\Theta_1\vcentcolon=\big\{\bs\te:\bs\te\in[-0.5+\epsilon,0.5]^q\big\}$ and $\Theta_2=\Omega_\beta\!\setminus\!\Theta_1$ (possibly an empty set), where $\Omega_\beta$ is given by \eqref{set}. Following Robinson (1995b) and Shimotsu (2007), we have
\small\begin{align}\label{diffp}
\P\big(\|\d-\do\|_\infty>\delta\big)&\leq\P\Big(\inf_{\overline{N_\delta}\cap\Omega_\beta}\!\!\big\{L(\bs d)\big\}\leq 0\Big)\nonumber\\
&\leq\P\Big(\inf_{\overline{N_\delta}\cap\Theta_1}\big\{L(\bs d)\big\}\leq 0\Big)+\P\Big(\inf_{\Theta_2}\big\{L(\bs d)\big\}\leq 0\Big)\vcentcolon=P_1+P_2,
\end{align}\normalsize
where, for a given set $\mathcal{O}$, $\overline{\mathcal O}$ denotes the closure of $\mathcal O$. We shall show that $P_1$ and $P_2$ go to zero as $n$ tends to infinity.  We deal with $P_1$ first. Notice that $L(\bs d)$ can be rewritten as
\small\begin{align}\label{aux1}
L(\bs d )&=\log\big(\det \{\widehat G(\bs d)\}\big)-\log\big(\det \{\widehat G(\do)\}\big)-2\sum_{k=1}^q\te_k\frac{1}{m}\sum_{j=1}^m\log(\l_j)\nonumber\\
&=\log\big(\det \{\widehat G(\bs d)\}\big)-\log\big(\det \{\widehat G(\do)\}\big)+\log\left(\frac{2\pi m}{n}\right)^{\!\!-2\sum_{k}\te_k}-\nonumber\\
&\hspace{1cm}-2\sum_{k=1}^q\te_k\bigg(\frac{1}{m}\sum_{j=1}^m\log(j)-\log(m)\bigg)-\sum_{k=1}^q\log(2\te_k+1)\nonumber\\
&=\log\big(\mz{A}(\bs d)\big)-\log\big(\mz{B}(\bs d)\big)-\log\big(\mz{A}(\do)\big)+\log\big(\mz{B}(\do)\big)+\mz{R}(\bs d)\nonumber\\
&= Q_1(\bs d)-Q_2(\bs d)+\mz{R}(\bs d),
\end{align}\normalsize
where
\small\begin{gather*}
Q_1(\bs d)\vcentcolon=\log\big(\mz{A}(\bs d)\big)-\log\big(\mz{B}(\bs d)\big),\qquad Q_2(\bs d)\vcentcolon=\log\big(\mz{A}(\do)\big)+\log\big(\mz{B}(\do)\big),\\
\mz{A}(\bs d)\vcentcolon =\left(\frac{2\pi m}{n}\right)^{\!\!-2\sum_{k}\te_k}\!\!\det \{\widehat G(\bs d)\},\qquad
\mz{B}(\bs d)\vcentcolon =\det \{G_0\}\prod_{k=1}^q\frac{1}{2\te_k+1},\\
\mbox{ and }\quad\mz{R}(\bs d)\vcentcolon =2\sum_{k=1}^q\te_k\bigg(\log(m)-\frac{1}{m}\sum_{j=1}^m\log(j)\bigg)-\sum_{k=1}^q\log(2\te_k+1).\phantom{ and and a}
\end{gather*}
\normalsize

By lemma 2 in Robinson (1995b), $\log(m)-m^{-1}\sum_{j=1}^m\log(j)=1+O(m^{-1}\log(m))$, so that
\small\[\mz{R}(\bs d)=\sum_{k=1}^q\Big[2\te_k-\log(2\te_k+1)\Big]+O\left(\frac{\log(m)}{m}\right).\]\normalsize
Since $x-\log(x+1)$ has a unique global minimum in $(-1,\infty)$ at $x=0$ and $x-\log(x+1)\geq x^2/4$, for $|x|\leq1$, it follows that
\small\begin{equation}\label{auxl1}\inf_{\overline{N_\delta}\cap\Theta_1}\!\!\big\{\mz{R}(\bs d)\big\}\geq \frac{1}{4}\Big(2\max_{k}\{\te_k\}\Big)^2\geq\delta^2>0\,.\end{equation}\normalsize

In view of \eqref{aux1} and \eqref{auxl1}, in order to show that $P_1\rightarrow0$ it suffices to show the existence of a function $h(\bs d)>0$ satisfying
\small\begin{equation}\label{lcond} (\mathrm{i})\,\,\sup_{\Theta_1}\big\{\big|\mz{A}(\bs d)-h(\bs d)\big|\big\}=o_{\P}(1);\quad(\mathrm{ii})\,\,h(\bs d)\geq \mz{B}(\bs d);\quad (\mathrm{iii})\,\,h(\do)=\mz{B}(\do), \end{equation}\normalsize
as $n$ goes to infinity, because (ii) implies $\displaystyle{\inf_{\Theta_1}}\big\{h(\bs d)\big\}\geq \displaystyle{\inf_{\Theta_1}}\big\{\mz{B}(\bs d)\big\}>0$, so that,
\small\[Q_1(\bs d)\geq \log\big(\mz{A}(\bs d)\big)-\log\big(h(\bs d)\big)=\log\big(h(\bs d)+o_\P(1)\big)-\log\big(h(\bs d)\big)=o_\P(1),\]\normalsize
and (iii) implies $Q_2(\bs d)=\log\big(h(\do)+o_\P(1)\big)-\log\big(h(\do)\big)=o_\P(1)$, both asymptotic orders being uniform in $\Theta_1$ in view of (i), and these results together with \eqref{auxl1} imply $P_1 \! \! \longrightarrow \! 0$.

To show (i), recall that
\small\begin{align}\label{lambda}
\L_j(\bs d)^{-1}&=\!\!\diag_{k\in\{1,\cdots,q\}}\!\!\{\l_j^{d_k}\mathrm{e}^{\im(\l_j-\pi)d_k/2}\}=\!\!
\diag_{k\in\{1,\cdots,q\}}\!\!\{\l_j^{(d_k-d_k^0)}\mathrm{e}^{\im(\l_j-\pi)(d_k-d_k^0)/2}\!\times\!\l_j^{d_k^0}\mathrm{e}^{\im(\l_j-\pi)d_k^0/2}\}\nonumber\\
&=\L_j(\bs d-\do)^{-1}\L_j(\do)^{-1} =\L_j(\bs \te)^{-1}\L_j(\do)^{-1}\!,
\end{align}\normalsize
so that we can write
\small\begin{align}\label{Mjs}
\mz A&(\bs d)=\left(\frac{2\pi m}{n}\right)^{\!\!-2\sum_{k}\te_k}\!\!\times\det \bigg\{\frac{1}{m}\sum_{j=1}^m\re\big[\L_j(\bs\te)^{-1}\L_j(\do)^{-1}f_n(\l_j)\comp{\L_j(\do)^{-1}} \comp{\L_j(\bs\te)^{-1}}\big]\bigg\}\nonumber\\
&=\det \bigg\{\frac{1}{m}\sum_{j=1}^m\re\big[M_j(\bs\te)\L_j(\do)^{-1}\big(f(\l_j)+o_\P(n^{-\beta})\big)\comp{\L_j(\do)^{-1}} \comp{M_j(\bs\te)}\big]\bigg\}\nonumber\\
&=\det \bigg\{\frac{1}{m}\sum_{j=1}^m\re\big[M_j(\bs\te)G_0\comp{M_j(\bs\te)}\big] +\frac{1}{m}\sum_{j=1}^m\re\big[M_j(\bs\te)\L_j(\do)^{-1}o_\P(n^{-\beta})\comp{\L_j(\do)^{-1}}\comp{M_j(\bs\te)}\big]\bigg\},
\end{align}\normalsize
where $M_j(\bs\te):=\diag_{k\in\{1,\cdots,q\}}\Big\{\mathrm{e}^{\im(\l_j-\pi){d}_k^0/2}(j/m)^{\te_k}\Big\}$. By omitting  the $\re[\cdot]$ operator, the $(r,r)$-th element in the second term inside the determinant on the right-hand side of \eqref{Mjs} is
\small\begin{align}\label{auxl2}
\frac 1m\sum_{j=1}^m\bigg(\frac jm\bigg)^{2\te_r}\l_j^{2d_r^0}o_\P(n^{-\beta})=o_\P(1)\frac 1m\sum_{j=1}^m\bigg(\frac jm\bigg)^{2\te_r}j^{2d_r^0}=o_\P(1)\frac{1}{m^{2d_r^0+1}}\sum_{j=1}^m\bigg(\frac jm\bigg)^{2d_r},
\end{align}\normalsize
since $\do\in\Omega_\beta$. Now, by taking $\gamma=2d_r+1>0$ in lemma 1 in Robinson (1995b), it follows that
\small\begin{align}\label{auxl3}
\frac 1m \sum_{j=1}^m\bigg(\frac jm\bigg)^{2d_r}&=\frac1{2d_r+1}\bigg[\frac{2d_r+1}{m}\sum_{j=1}^m\bigg(\frac jm\bigg)^{2d_r}\bigg]=\frac1{2d_r+1}\big[O(m^{\beta-1})+1\big]
\end{align}\normalsize
so that \eqref{auxl2} is $o_\P(1)$ uniformly in $\Theta_1$ in view of \eqref{auxl3} and since $\do\in\Omega_\beta$ by hypothesis and $\bs d\in \Omega_\beta$ by the definition of the estimator. Hence, the second term inside the determinant on the right-hand side of \eqref{Mjs} is $o_\P(1)$ uniformly in $\Theta_1$, so that
\small\begin{align}\label{auxl4}
\mz A&(\bs d)=\det \bigg\{\frac{1}{m}\sum_{j=1}^m\re\big[M_j(\bs\te)G_0\comp{M_j(\bs\te)}\big] +o_\P(1)\bigg\}.
\end{align}\normalsize
Upon defining the matrices
\small\[\mathcal{E}_0:=\diag_{k\in\{1,\cdots,q\}}\Big\{\mathrm e^{\im\pi d_k^0/2}\Big\}, \quad\mathcal G_0:=\re\big[\mathcal{E}_0G_0\comp{\mathcal E_0}\big]\ \ \mbox{ and }\ \ \mathcal{M}(\bs\te):=\left(\frac{1}{1+\te_r+\te_s}\right)_{r,s=1}^q\]\normalsize
from the proof of theorem 3(a) in Shimotsu (2007), it follows that the function
\small\[h(\bs d):=\det \Big\{\mathcal M(\bs\te)\odot \mathcal G_0\Big\},\]\normalsize
where $\odot$ denotes the Hadamard product, satisfies the conditions (i), (ii) and (iii) in \eqref{lcond} (notice that \eqref{auxl4} is the same equation as the one following (40) in Shimotsu, 2007, p.303, with the obvious identifications).

Now we move to bound $P_2$. By using \eqref{lambda}, rewrite $L(\bs d)$ as
\small\begin{align}\label{aux6}
L(\bs d )=\log\Big(\det \big\{\mcd(\bs d)\big\}\Big)-\log\Big(\det \big\{\mcd(\do)\big\}\Big),
\end{align}\normalsize
where
\small\[\mcd(\bs d)\vcentcolon=\frac{1}{m}\sum_{j=1}^m\re\big[\mcp_j(\bs \te)\L_j(\do)^{-1}f_n(\l_j)\comp{\L_j(\do)^{-1}}\comp{\mcp_j(\bs \te)}\big],\]\normalsize
with
\small\[  \mcp_j(\bs\te)\vcentcolon=\diag_{k\in\{1,\cdots,q\}}\bigg\{\mathrm{e}^{\im (\l_j-\pi)d^0_k/2}\left(\frac{j}{\mz p}\right)^{\!\!\te_k}\bigg\}\quad\mbox{and}\quad\mz p\vcentcolon=\exp\bigg(\frac{1}{m}\sum_{j=1}^m\log(j)\bigg),\]\normalsize
and, as $m$ tends to infinity, $\mz p\sim m/\mathrm{e}$. Observe that $\mcd(\bs d)$ is positive semidefinite since each summand of $\mcd$ is. For $\kappa\in (0,1)$, define the auxiliary functions
\small\[\mdk(\bs d)\vcentcolon=\sz\!\re\big[\mcp_j(\bs \te)\L_j(\do)^{-1}f_n(\l_j)\comp{\L_j(\do)^{-1}}\comp{\mcp_j(\bs \te)}\big]\]\normalsize
and
\small\[ \qz\vcentcolon=\sz\!\re\big[\mcp_j(\bs\te)G_0\comp{\mcp_j(\bs\te)}\big],\]\normalsize
where $\lfloor x\rfloor$ denotes the integer part of  $x$. Now, by the hypothesis on $f_n$,
\small\begin{align}\label{auxl5}
\mdk(\bs d)&=\sz\!\re\big[\mcp_j(\bs \te)\L_j(\do)^{-1}\big(f(\l_j)+o_\P(n^{-\beta})\big)\comp{\L_j(\do)^{-1}}\comp{\mcp_j(\bs \te)}\big]\nonumber\\
&=\qz+o(1)+\sz\re\big[\mcp_j(\bs \te)\L_j(\do)^{-1}o_\P(n^{-\beta})\comp{\L_j(\do)^{-1}}\comp{\mcp_j(\bs \te)}\big],
\end{align}\normalsize
where the last equality follows from lemma 5.4 in Shimotsu and Phillips (2005).
The $(r,s)$-th element of the third term on the RHS of \eqref{auxl5} is given by
\small\begin{align*}
\re\bigg[\sz  \left(\frac{j}{\mz p}\right)^{\!\!\te_r+\te_s}&\left(\frac{2\pi j}{n}\right)^{\!\!d^0_r+d^0_s}\mathrm{e}^{\im (\l_j-\pi)(\te_r-\te_s)/2}o_\P(n^{-\beta})\bigg]=\\
&=O(1)\left(\frac{m}{\mz p}\right)^{\!\!\te_r+\te_s}\left(\frac{m}{n}\right)^{\!\!d^0_r+d^0_s}o_\P(n^{-\beta})\,\sz \left(\frac{j}{m}\right)^{\!\!2(d^0_r+d^0_s)-(\widehat d_r+\widehat d_s)}\\
&=O(1)o_\P(1)\left(\frac{m}{\mz p}\right)^{\!\!\te_r+\te_s}O(1)=o_\P(1),
\end{align*}\normalsize
uniformly in $\bs \te\in\Theta_2$, since $\do\in\Omega_\beta$, $\beta\in(0,1)$ and Assumption \textbf{A4}, where the penultimate equality follows from lemma 5.4 in Shimotsu and Phillips (2005). This shows that
\small\[\sup_{\Theta_2}\big\{\big|\det\{\mcd(\bs d)\}-\det\{\qz\}\big|\big\}=o_\P(1).\]\normalsize
The proof now follows viz a viz (with the obvious notational identification) from the proof of theorem 3(a) in Shimotsu (2007), p.303 (see the argument following equation (42)). From this, we conclude that $P_2\longrightarrow0$, as $n$ tends to infinity, and this completes the proof. \fim

\subsubsection{Proof of Lemma \ref{lemma2}}

For fixed $r,s\in\{1,\cdots,q\}$, let $\mathscr{A}_{uv}\vcentcolon=\sum_{j=u}^v\mathscr A_j$ and $\mathscr{B}_{uv}\vcentcolon=\sum_{j=u}^v\mathscr B_j$, where
\small\begin{equation}\label{aux8}
\mathscr A_j\vcentcolon=\mathrm{e}^{\im(\l_j-\pi)(d_r^0-d_s^0)/2}\l_j^{d^0_r+d^0_s}\big[f_n^{rs}(\l_j) - \big(A(\l_j)\big)_{r\bs\cdot}I_{\bs\eps}(\l_j) \big(\comp{A(\l_j)}\big)_{\bs\cdot s}\big],
\end{equation}\normalsize
and
\small\begin{equation}\label{aux9}
\mathscr{B}_j\vcentcolon=\mathrm{e}^{\im(\l_j-\pi)(d_r^0-d_s^0)/2}\l_j^{d^0_r+d^0_s}\big(A(\l_j)\big)_{r\bs\cdot}I_{\bs\eps}(\l_j) \big(\comp{A(\l_j)}\big)_{\bs\cdot s}-G_0^{rs}.
\end{equation}\normalsize
Hence, for each $j$, $\mathscr A_j+\mathscr B_j=\mathrm{e}^{\im(\l_j-\pi)(d_r^0-d_s^0)/2}\l_j^{d^0_r+d^0_s}f_n^{rs}(\l_j)-G_0^{rs}$. For fixed $u\leq j\leq v$, we have
\small\begin{align*}
\E\big(|\mathscr A_j|\big)&=(2\pi)^{d_r^0+d_s^0}\E\Big(\l_j^{-d^0_r-d^0_s}\Big|f_n^{rs}(\l_j) - \big(A(\l_j)\big)_{r\bs\cdot}I_{\bs\eps}(\l_j) \big(\comp{A(\l_j)}\big)_{\bs\cdot s}\Big|\Big)\ =\ o(1)\\
&=(2\pi )^{d_r^0+d_s^0}o(1)\,\,=\,\,o(1),
\end{align*}\normalsize
from where we conclude that $\displaystyle{\max_{r,s}}$ $\!\big\{\E\big(\big|\sum_{j=u}^v\mathscr{A}_j\big|\big)\big\} = o(v-u+1)$ and the result on $\mathscr A_{uv}$ follows.
As for $\mathscr B_j$, from the proof of lemma 1(a) in Shimotsu (2007) (notice that $\mathscr B_j$ does not depend on $f_n$) it follows that $\sum_{j=u}^v\mathscr B_j=o_\P(v)$ uniformly in $u$ and $v$ and hence the desired result on $\mathscr{B}_{uv}$ follows. The last assertion is now straightforward. \fim

\subsubsection{Proof of Theorem \ref{nonbeta}}

By carefully inspecting the proof of Theorem \ref{consistency}, we observe that it suffices to show (with the same notation as in the aforementioned proof) \eqref{auxl4} uniformly in $\Theta_1$ and that $\mdk(\bs d)-\qz=o_\P(1)$ uniformly in $\Theta_2$. To show \eqref{auxl4},  it suffices to show that
\small\begin{equation}\label{mn1}
\frac1m\sum_{j=1}^m\re\Big[M_j(\bs\te)\Lambda_j(\do)^{-1}f_n(\l_j)\comp{\Lambda_j(\do)^{-1}}\comp{M_j(\bs\te)}\Big] =\frac1m\sum_{j=1}^m\re\Big[M_j(\bs\te)G_0\comp{M_j(\bs\te)}\Big]+ o_\P(1).
\end{equation}\normalsize
The $(r,s)$-th component of the LHS of \eqref{mn1} is given by
\small\[\frac{1}{m}\sum_{j=1}^m\re\bigg[\mathrm{e}^{\im(\l_j-\pi)(d^0_r-d^0_s)/2}\left(\frac{j}{m}\right)^{\te_r+\te_s} \!\! f_n^{rs}(\l_j)\Big(\L_j^{(r)}(\do)\comp{\L_j^{(s)}(\do)}\Big)^{-1}\bigg].\]\normalsize
Summation by parts (see Zygmund, 2002, p.3) yields
\small
\begin{align}\label{op1}
\sup_{\Theta_1}\bigg\{\bigg|\frac{1}{m}&\sum_{j=1}^m\mathrm{e}^{\im(\l_j-\pi)(d^0_r-d^0_s)/2}\left(\frac{j}{m}\right)^{\te_r+\te_s} \!\! \Big[f_n^{rs}(\l_j)\Big(\L_j^{(r)}(\do)\comp{\L_j^{(s)}(\do)}\Big)^{\!\!-1}\!\!\!-G_0^{rs}\Big]\bigg|\bigg\}\leq\nonumber\\
&\leq\frac{1}{m}\sum_{k=1}^{m-1}\sup_{\Theta_1}\bigg\{\bigg|\mathrm{e}^{\im(\l_k-\pi)(d^0_r-d^0_s)/2} \left(\frac{k}{m}\right)^{\te_r+\te_s}\hspace{-.3cm}-\mathrm{e}^{\im(\l_{k+1}-\pi)(d^0_r-d^0_s)/2} \left(\frac{k+1}{m}\right)^{\te_r+\te_s}\bigg|\bigg\}\times\nonumber\\
&\hspace{2cm}\times\,\,\,\bigg|\sum_{j=1}^k\Big[f_n^{rs}(\l_j)\Big(\L_j^{(r)}(\do)\comp{\L_j^{(s)}(\do)}\Big)^{\!\!-1}\!\!\!-G_0^{rs}\Big]\bigg|\ +\nonumber\\
&\hspace{3cm}+\,\,\,\bigg|\frac{1}{m}\sum_{j=1}^m\Big[f_n^{rs}(\l_j)\Big(\L_j^{(r)}(\do)\comp{\L_j^{(s)}(\do)}\Big)^{\!\!-1}\!\!\!-G_0^{rs}\Big]\bigg|\nonumber\\
&\leq\mz C\sum_{k=1}^{m-1}\left(\frac{k}{m}\right)^{\!\!2\epsilon}\frac{1}{k^2}\bigg|\sum_{j=1}^k\Big[f_n^{rs}(\l_j)\Big(\L_j^{(r)}(\do) \comp{\L_j^{(s)}(\do)}\Big)^{\!\!-1}\!\!\!-G_0^{rs}\Big]\bigg|\ +\nonumber\\
&\hspace{2cm}+\,\,\,\bigg|\frac{1}{m}\sum_{j=1}^m\Big[f_n^{rs}(\l_j)\Big(\L_j^{(r)}(\do)\comp{\L_j^{(s)}(\do)}\Big)^{\!\!-1}\!\!\!-G_0^{rs}\Big]\bigg|,
\end{align}
\normalsize
where $0<\mz C<\infty$ is a constant. The first term in \eqref{op1} is, by \eqref{condthm},
\small\begin{align}\label{new1}
\sum_{k=1}^{m-1}\left(\frac{k}{m}\right)^{2\epsilon}\!\!\!\frac{1}{k^2}\bigg|\sum_{j=1}^k\Big[f_n^{rs}(\l_j)&\Big(\L_j^{(r)}(\do) \comp{\L_j^{(s)}(\do)}\Big)^{\!\!-1}\!\!\!-G_0^{rs}\Big]\bigg|\leq\sum_{k=1}^{m-1}\left(\frac{k}{m}\right)^{2\epsilon}\!\!\!\frac{1}{k^2}\Big(\big|\mathscr{A}_{1k}\big| +\big|\mathscr{B}_{1k}\big|\Big)\nonumber\\
&=\frac{1}{m^{2\epsilon}}\sum_{k=1}^{m-1}k^{2(\epsilon-1)}\big|\mathscr{A}_{1k}\big|+\frac{1}{m^{2\epsilon}}\sum_{k=1}^{m-1}k^{2(\epsilon-1)}o_\P(k) = o_\P(1),
\end{align}\normalsize
uniformly in $(r,s)$, where the last equality follows from lemma 2 in Robinson (1995b), Lemma \ref{lemma2} and Chebyshev's inequality, since
\small\begin{align*}
\E\left(\frac{1}{m^{2\epsilon}}\sum_{k=1}^{m-1}k^{2(\epsilon-1)}\big|\mathscr{A}_{1k}\big|\right)&=\frac{2\epsilon}{m}\, \sum_{k=1}^{m-1}{\left(\frac km\right)}^{2\epsilon-1}\frac{1}{k}\, \E\Big(\big|\mathscr{A}_{1k}\big|\Big)=\big[1+O(m^{-\epsilon})\big]o(1)= o(1).
\end{align*}\normalsize
%
%
The second term in \eqref{op1} is $o_\P(1)$ uniformly in $(r,s)$ by \eqref{condthm}. Hence, \eqref{mn1} follows. As for the difference $\mdk(\bs d)-\qz$, its $(r,s)$-th element is given by
\small
\begin{align}
\sz\!&\re\bigg[\mathrm{e}^{\im(\l_j-\pi)(d^0_r-d^0_s)/2}\left(\frac{j}{\mz p}\right)^{\te_r+\te_s} \!\! \Big[f_n^{rs}(\l_j)\Big(\L_j^{(r)}(\do)\comp{\L_j^{(s)}(\do)}\Big)^{\!\!-1}\!\!\!-G_0^{rs}\Big]\bigg]=\nonumber\\
&=\left(\frac{m}{\mz p}\right)^{d^0_r+d^0_s} \!\!\re\bigg[\sz\mathrm{e}^{\im(\l_j-\pi)(\te_r-\te_s)/2}\left(\frac{j}{m}\right)^{\te_r+\te_s} \!\!\!\times\Big[f_n^{rs}(\l_j)\Big(\L_j^{(r)}(\do)\comp{\L_j^{(s)}(\do)}\Big)^{\!\!-1}\!\!\!-G_0^{rs}\Big]\bigg]\nonumber
\end{align}
\normalsize
which is $o_\P(1)$ uniformly in $\bs\te\in\Theta_2$ by similar argument as the one applied in deriving \eqref{op1}, from summation by parts and lemma 5.4 in Shimotsu and Phillips (2005). This completes the proof.\fim

\subsubsection{Proof of Lemma \ref{lemma1b2}}

\rm{(a)} For $r,s\in\{1,\cdots,q\}$ fixed, ignoring the maximum in expression \eqref{lemma1b2a} for a while, we see that if $d_r^0+d_s^0\geq0$,
\small\begin{align*}
\sum_{j=1}^v\l_j^{d_r^0+d_s^0}&\Big[f_n^{rs}(\l_j)- \big(A(\l_j)\big)_{r\bs\cdot}I_{\bs\eps}(\l_j)\big(\comp{A(\l_j)}\big)_{\bs\cdot s}\Big]\leq\\
&\leq\ \max\big\{1,n^{|d_r^0+d_s^0|}\big\}\max_{1\leq v\leq m}\bigg\{\sum_{j=1}^v\Big[f_n^{rs}(\l_j)- \big(A(\l_j)\big)_{r\bs\cdot}I_{\bs\eps}(\l_j)\big(\comp{A(\l_j)}\big)_{\bs\cdot s}\Big]\bigg\}\\
&=\ o_\P\left(\frac {m^2}{n}\right)\ =\ o_\P\left(\frac{\sqrt{m}}{\log(m)}\right),  
\end{align*}\normalsize
where the last equality follows since Assumption \textbf{B4} holds for $\alpha=1$, which implies $mn^{-1}=o\big(m^{-1/2}\log(m)^{-1}\big)$.

\noindent\rm{(b)}  Rewrite the argument of the summation in \eqref{lemma1b2b} as $\mathscr{A}_j+\mathscr{B}_j+\mathscr{C}_j$, where
\small\begin{align*}
\mathscr{A}_j&:=\l_j^{d_r^0+d_s^0}\Big[f_n^{rs}(\l_j)- \big(A(\l_j)\big)_{r\bs\cdot}I_{\bs\eps}(\l_j)\big(\comp{A(\l_j)}\big)_{\bs\cdot s}\Big],\\
\mathscr{B}_j&:=\l_j^{d_r^0+d_s^0}\Big[ \big(A(\l_j)\big)_{r\bs\cdot}I_{\bs\eps}(\l_j)\big(\comp{A(\l_j)}\big)_{\bs\cdot s}-f_{rs}(\l_j)\Big],\\
\mathscr{C}_j&:=\l_j^{d_r^0+d_s^0}f_{rs}(\l_j)-\mathrm{e}^{\im(\pi-\l_j)(d_r^0-d_s^0)/2}G_0^{rs}.
\end{align*}\normalsize
Part \rm{(a)} yields $\max_{v}\Big\{\max_{r,s}\Big\{\!\!\sum_{j=1}^v\big|\mathscr{A}_j\big|\Big\}\Big\} =o_\P\Big(m^{1/2}\big(\log(m)\big)^{-1}\Big)$, while, from the proof of lemma 1(b2) in Shimotsu (2007), we obtain $\max_{r,s}\Big\{\!\!\sum_{j=1}^v\big|\mathscr{B}_j\big|\Big\} =O_\P\big(v^{1/2}\log(v)\big)$ (recall that $\mathrm{e}^{\im\l_j}=O(\l_j)$). Assumption \textbf{B1} implies that $|\mathscr C_j|=O(\l_j^{\alpha})$ regardless $r,s$ so that $\max_{r,s}\Big\{\!\!\sum_{j=1}^v\big|\mathscr{C}_j\big|\Big\} =O\big(v^{\alpha+1}n^{-\alpha}\big)$. The result now follows by noticing that $m^{1/2}\big(\log(m)\big)^{-1}=O\big(m^{1/2}\log(m)\big)$.\fim

\subsubsection{Proof of Theorem \ref{norm}}

The argument is similar to the one in the proof of theorem 3(b) in Shimotsu (2007). By hypothesis, with probability tending to 1 as $n$ tends to infinity,
\small\[\bs0=\frac{\partial S(\bs d)}{\partial\bs d}\bigg|_{\d}=\frac{\partial S(\bs d)}{\partial\bs d}\bigg|_{\do}+\bigg(\frac{\partial^2S(\bs d)}{\partial \bs d\partial\bs{d}^{\,\prime}}\bigg|_{\od}\bigg)(\d-\do),\]\normalsize
where $\od\in\R^q$ is such that $\|\od-\do\|_{\infty}\leq\|\d-\do\|_{\infty}$. Notice that $\d$ has the stated limiting distribution if
\small\begin{equation}\label{M1}
\sqrt m\ \frac{\partial S(\bs d)}{\partial\bs d}\bigg|_{\do}\overset{d}{-\!\!\!\longrightarrow} N(0,\Sigma)
\end{equation}\normalsize
and
\small\begin{equation}\label{M2}
\frac{\partial^2S(\bs d)}{\partial \bs d\partial\bs{d}^{\,\prime}}\bigg|_{\od}\overset{\P}{-\!\!\!\longrightarrow}\Omega.
\end{equation}\normalsize
We start by proving \eqref{M1}. First notice that
\small\begin{align*}
\sqrt{m}\ \frac{\partial S(\bs d)}{\partial\bs d}\bigg|_{\do}=-\frac{2}{\sqrt{m}}\sum_{j=1}^m \log(\l_j)+\tr\bigg[\G{\do}^{-1}\sqrt{m}\,\frac{\partial\G{\do}}{\partial d_r}\bigg].
\end{align*}\normalsize
Let $\ir $ denote a $q\times q$ matrix whose $(r,r)$-th element is 1 and all other elements are zero and let $\mathcal E_j:=\displaystyle{\diag_{k=1,\cdots,q}\{\mathrm e^{\im(\pi-\l_j)d_k^0/2}\}}$. With this notation, for $r\in\{1,\cdots,q\}$ we can write
\small\[\G{\do}=\frac1m\sum_{j=1}^m\re\Big[\mathcal E_j\L_j(\do)^{-1}f_n(\l_j)\comp{\L_j(\do)^{-1}}\comp{\mathcal E_j}\Big]\]\normalsize
and
\small\[\mathscr H(r):=\sqrt{m}\ \frac{\partial\G{\bs d}}{\partial d_r}\bigg|_{\do}=\frac{1}{\sqrt{m}}\sum_{j=1}^m\log(\l_j)\re\Big[\mathcal E_j\L_j(\do)^{-1}\big(\ir f_n(\l_j)+f_n(\l_j)\ir \big) \comp{\L_j(\do)^{-1}}\comp{\mathcal E_j}\Big].\]\normalsize
For $j\in\{1,\cdots,m\}$, let
\small\[a_j\vcentcolon=\log(\l_j)-\frac{1}{m}\sum_{k=1}^m\log(\l_k)=\log(j)-\frac{1}{m}\sum_{k=1}^m\log(k)=O\big(\log(m)\big),\]\normalsize
and recall that for a matrix $M$, $(M)_{r\bs\cdot}$ and $(M)_{\bs\cdot s}$  denote, respectively, the $r$-th row and the $s$-th column of $M$. To show \eqref{M1}, we shall apply the Cr\'amer-Rao device. For an arbitrary vector $\bs\eta\in\R^q$, we have
\[{\bs\eta}^\prime\sqrt m\ \frac{\partial S(\bs d)}{\partial \bs d}\bigg|_{\do}=\sum_{k=1}^q\eta_k\bigg[\sqrt m\ \frac{\partial S(\bs d)}{\partial d_k}\bigg|_{\do}\bigg].\]
Observe that we can write
\begin{align}\label{n3}
\sqrt m\ \frac{\partial S(\bs d)}{\partial d_k}\bigg|_{\do}&=\tr\bigg[\G{\do}^{-1}\bigg(\mathscr H(k)-\frac{2}{\sqrt{m}}\sum_{j=1}^m\log(\l_j)\G{\do}\mathrm{I}_{(k)}\bigg)\bigg]\nonumber\\
&=\tr\bigg[\G{\do}^{-1}\frac{2}{\sqrt{m}}\sum_{j=1}^ma_j\re\Big[\lf\Big]\mathrm{I}_{(k)}\bigg]\nonumber\\
&=\Big[(G_0^{-1})_{k\bs\cdot}+o_\P(1)\Big]\frac{2}{\sqrt{m}}\sum_{j=1}^ma_j\Big(\re\Big[\lf\Big]\Big)_{\bs\cdot k},
\end{align}
where the last equality follows from Lemma 3.1(b). Omitting the $\re[\cdot]$ operator and upon rewriting $f_n(\l_j)=f_n(\l_j)-A(\l_j)I_{\bs\eps}(\l_j)\comp{A(\l_j)}+A(\l_j)I_{\bs\eps}(\l_j)\comp{A(\l_j)}$, the summation in \eqref{n3} can be split into two parts, one of them reads
\small
\begin{align*}
\bigg|\bigg(\sum_{j=1}^ma_j\mathcal E_j\L_j(\do)^{-1}\Big(&f_n(\l_j)-A(\l_j)I_{\bs\eps}(\l_j)\comp{A(\l_j)}\Big)\comp{\L_j(\do)^{-1}}\comp{\mathcal E_j}\bigg)_{rs}\bigg|\leq\\
&\leq O\big(\log(m)\big)\max_{v=1,\cdots,m}\!\!\bigg\{\Big|\sum_{j=1}^v\l_j^{d_r^0+d_s^0}
\Big(f_n^{rs}(\l_j)-\big(A(\l_j)\big)_{r\bs\cdot}I_{\bs\eps}(\l_j)\big(\comp{A(\l_j)}\big)_{\bs\cdot s}\Big)\Big|\bigg\}\\
&\hskip-.2cm=O\big(\log(m)\big)o_\P\bigg(\frac{\sqrt{m}}{\log(m)}\bigg)\,\,=\,\,o_\P(\sqrt{m}),
\end{align*}
\normalsize
uniformly in $r,s\in\{1,\cdots,q\}$ by Lemma 3.1(a), so that
\small\begin{align}\label{n4}
\frac{1}{\sqrt{m}}\sum_{j=1}^ma_j&\lf=\nonumber\\
&=\frac{1}{\sqrt{m}}\sum_{j=1}^ma_j \mathcal E_j\L_j(\do)^{-1}A(\l_j)I_{\bs\eps}(\l_j)\comp{A(\l_j)}\comp{\L_j(\do)^{-1}}\comp{\mathcal E_j} + \frac{1}{\sqrt{m}}\,o_\P(\sqrt{m})\nonumber\\
&=\frac{1}{\sqrt{m}}\sum_{j=1}^ma_j \Big[\mathcal E_j\L_j(\do)^{-1}A(\l_j)I_{\bs\eps}(\l_j)\comp{A(\l_j)}\comp{\L_j(\do)^{-1}}\comp{\mathcal E_j}-\mathcal E_jG_0\comp{\mathcal E_j}\Big]+o_\P(1),
\end{align}\normalsize
where the last equality follows from $\sum_{j=1}^m a_j=0$. Notice that \eqref{n4} no longer depends on $f_n$ and, with the obvious notational identification, equation \eqref{n4} is exactly equation (45) on the proof of theorem 3(b), subsection A.4.1 in Shimotsu 2007, p.305. The proof of \eqref{M1} now follows viz a viz from the proof of that theorem.

We proceed to show \eqref{M2}. Let $\bs\te:=\bs d-\do$ and define
\small\[\mathcal M:=\Big\{\bs d:\log(n)^4\|\bs d -\do\|_{\infty}<\delta\Big\}=\Big\{\bs \te :\log(n)^4\|\bs\te\|_{\infty}<\delta\Big\}.\]\normalsize
First we show that $\P\big(\od\in \mathcal M\big)\rightarrow 1$, as $n\rightarrow\infty$.   Assuming the same notation as in Theorem \ref{consistency}, recall the decomposition of $L(\bs d)=S(\bs d)-S(\do)= Q_1(\bs d)-Q_2(\bs d)+\mz{R}(\bs d)$ given in \eqref{aux1}. The same argument as in the proof of Theorem \ref{consistency} yields
\small\[\inf_{\Theta_1\setminus \mathcal M}\big\{\mz R(\bs d)\big\}\geq\delta^2\log(n)^8,\]\normalsize
and upon applying  Lemma 3.1, we obtain
\small\[\sup_{\Theta_1}\Big\{\big|\mz A(\bs d)-h(\bs d)\big|\Big\}=O_\P\bigg(\frac{m^\alpha}{n^\alpha}+\frac{\log(m)}{m^{2\epsilon}}+\frac{m}{n}\bigg).\]\normalsize
It follows, uniformly in $\Theta_1$ (cf. Shimotsu, 2007, p.300), that
\small\begin{align*}
\log\big(\mz A(\bs d)\big)&-\log\big(\mz B(\bs d)\big)\geq \log\big(h(\bs d)+o_\P\big(\log(n)^{-8}\big)\big)-\log\big(h(\bs d)\big)=o_\P\big(\log(n)^{-8}\big)\\
\log\big(\mz A(\do)\big)&-\log\big(\mz B(\do)\big)=\log\big(h(\do)+o_\P\big(\log(n)^{-8}\big)\big)-\log\big(h(\do)\big)=o_\P\big(\log(n)^{-8}\big),
\end{align*}\normalsize
from where we conclude that $\P\big(\inf_{\Theta_1\setminus\mathcal M}L(\bs d)\leq 0\big)\rightarrow0$. Hence $\P\big(\od\in \mathcal M\big)\rightarrow 1$, as $n\rightarrow\infty$. Next we observe that
\small\[\frac{\partial^2 S(\bs d)}{\partial d_r\partial d_s}=\tr\bigg[-\G{\bs d}^{-1}\,\frac{\partial\G{\bs d}}{\partial d_r}\,\G{\bs d}^{-1}\frac{\partial\G{\bs d}}{\partial d_s}+\G{\bs d}^{-1}\frac{\partial^2 \G{\bs d}}{\partial d_r\partial d_s}\bigg].\]\normalsize
For $k\in\{0,1,2\}$, let
\small\begin{align*}
J_k(\bs d)&:=\frac1m\sum_{j=1}^m\log(\l_j)^k\re\big[\lf\big].
\\&=\frac1m\sum_{j=1}^m\log(\l_j)^k\re\Big[\diag_{i=1,\cdots,q}\{\l_j^{d_i}\}f_n(\l_j)\diag_{i=1,\cdots,q}\{\l_j^{d_i}\}\Big].
\end{align*}\normalsize
We notice that $J_0(\bs d)=\G{\bs d}$ and the derivatives of $\G{\bs d}$ are given by
\small\begin{equation}\label{a35}
\frac{\partial\G{\bs d}}{\partial d_r}=\mathrm I_{(r)}J_1(\bs d)+J_1(\bs d)\mathrm I_{(r)}
\end{equation}\normalsize
and
\small\begin{equation}\label{a36}
\frac{\partial^2\G{\bs d}}{\partial d_r\partial d_s}=\mathrm I_{(r)}\mathrm I_{(s)}J_2(\bs d)+\mathrm I_{(r)}J_2(\bs d)\mathrm I_{(s)}+\mathrm I_{(s)}J_2(\bs d)\mathrm I_{(r)}+J_2(\bs d)\mathrm I_{(r)}\mathrm I_{(s)}.
\end{equation}\normalsize
From the proof of theorem 3(b) in Shimotsu (2007), it suffices to show that
\small\begin{equation}\label{fk}
J_k(\bs d)=\mathcal G_0\frac1m\sum_{j=1}^m\log(\l_j)^k+o_\P\big(\log(n)^{k-2}\big),
\end{equation}\normalsize
uniformly in $\bs d\in\mathcal M$ (notice that it also implies $\widehat G(\d)\overset{\P}{\longrightarrow}\mathcal G_0$). In order to do that, let
\small\[D_k(\bs\te):=\frac1m\sum_{j=1}^m\log(\l_k)^k\mathcal E_j\L_j(\bs\te)^{-1}G_0\comp{\L_j(\bs\te)^{-1}}\comp{\mathcal E_j}\]\normalsize
and notice that \eqref{fk} follows if
\small\begin{equation}\label{h5}
\sup_{\bs d\in \mathcal M}\bigg\{\bigg\|\frac{1}{m}\sum_{j=1}^m\log(\l_j)^k\mathcal E_j\L_j(\bs d)^{-1}f_n(\l_j)\comp{\L_j(\bs d)^{-1}}\comp{\mathcal E_j}-D_k(\bs \te)\bigg\|_\infty\bigg\} =o_\P\big(\log(n)^{k-2}\big),
\end{equation}\normalsize
and
\small\begin{equation}\label{h6}
\sup_{\bs d\in \mathcal M}\bigg\{\bigg\|D_k(\bs \te)-\mathcal G_0\frac{1}{m}\sum_{j=1}^m\log(\l_j)^k\bigg\|_\infty\bigg\} =o\big(\log(n)^{k-2}\big).
\end{equation}\normalsize
By applying \eqref{lambda}, \eqref{h5} can be rewritten as
\small\[\sup_{\bs d\in \mathcal M}\bigg\{\bigg\|\frac{1}{m}\sum_{j=1}^m\log(\l_j)^k\mathcal E_j\L_j(\bs \te)^{-1}\Big[\lfl-G_0\Big]\comp{\L_j(\bs \te)^{-1}}\comp{\mathcal E_j} \bigg\|_\infty\bigg\}.\]\normalsize
Define $b_j(\bs\te;k)\vcentcolon=\log(\l_j)^k\mathrm{e}^{\im(\l_j-\pi)(\te_r-\te_s)/2}\l_j^{\te_r+\te_s}$, for $k=0,1,2$. Then, by omitting the supremum, the $(r,s)$-th element of \eqref{h5} is equal to
\small\begin{align}\label{h7}
\bigg|\frac{1}{m}\sum_{j=1}^mb_j(\bs\te;k)&\Big[\l_j^{d^0_r+d^0_s}f_n^{rs}(\l_j)-\mathrm{e}^{\im(\pi-\l_j)(d_r^0-d_s^0)/2}G_0^{rs}\Big] \bigg|\leq\nonumber\\
&\leq\frac{1}{m}\sum_{j=1}^{m-1}\Big|b_j(\bs\te;k)-b_{j+1}(\bs\te;k) \Big| \bigg|\sum_{l=1}^{j}\l_l^{d^0_r+d^0_s} f_n^{rs}(\l_l)-\mathrm{e}^{\im(\pi-\l_l)(d_r^0-d_s^0)/2}G_0^{rs}\bigg|\nonumber\\
&\hspace{1cm}+\frac{b_m(\bs \te;k)}{m}\,\bigg|\sum_{j=1}^m\l_j^{d^0_r+d^0_s} f_n^{rs}(\l_j)-\mathrm{e}^{\im(\pi-\l_j)(d_r^0-d_s^0)/2}G_0^{rs}\bigg|,
\end{align}\normalsize
where the inequality follows from summation by parts. Now, since
\small\[b_j(\bs\te;k)-b_{j+1}(\bs\te;k) =O\left(\frac{\log(n)^k}{j}\right)\quad\mbox{ and }\quad b_m(\bs \te;k)=O\big(\log(n)^{k}\big),\]\normalsize
uniformly in $\bs \te\in\mathcal M$, for any $k=0,1,2$, by Lemma \ref{lemma1b2} and Remark \ref{remark} it follows that the first term on the RHS of \eqref{h7} is equivalent to
\small\begin{align*}
O\left(\frac{\log(n)^k}{m}\right)&\frac{1}{m}\,O_\P\left(\frac{m^{\alpha+1}}{n^\alpha}+m^{1/2}\log(m)\right)+\\
&\hspace{1cm}+ O\big(\log(n)^k\big)\frac{1}{m}\,O_\P\left(\frac{m^{\alpha+1}}{n^\alpha}+m^{1/2}\log(m)\right) =o_\P\big(\log(n)^{k-2}\big),
\end{align*}\normalsize
where the last equality follows from Assumption \textbf{B4} (see also Remark \ref{remark}), because
\small\begin{align*}
\log(n)^2\frac{1}{m}\,O_\P\left(\frac{m^{\alpha+1}}{n^\alpha}+m^{1/2}\log(m)\right)&=
\bigg[\frac{\log(n)^2}{m^{1/2}\log(m)}+\frac{\log(m)}{m^{1/4}}\,\frac{\log(n)^2}{m^{1/4}}\bigg]O_\P(1)=o_\P(1).
\end{align*}\normalsize
The other term in \eqref{h7} is dealt analogously, so that \eqref{h5} follows. As for \eqref{h6}, it does not depend on $f_n$ so that it follows from the proof of theorem 3(b) in Shimotsu (2007) (see section A.4.2, p.307) with the obvious notational adaptations. This completes the proof of \eqref{M2} and finishes the proof of the theorem.\fim

\subsubsection{Proof of Corollary \ref{tapered.cons}}
We shall show \eqref{indir}. From the proof of Lemma \ref{lemma2}, it suffices to show that, for $\mathscr A_j$ as in \eqref{aux8} and $f_n$ as in the enunciate, $\E(|\mathscr A_j|)=o(1)$ uniformly in $j$. From the proof of lemma 1(a) in Shimotsu (2007) (see also theorem 2 in Robinson, 1995a), we have
\small
\begin{subequations}
\begin{align}
&\E\big(I_n(\l_j)\big)=f(\l_j)\bigg(1+O\left(\frac{\log(j+1)}{j}\right)\bigg); \label{c1a}\\
&\E\big(I_{\bs\eps}(\l_j)\big)=\frac{\mathrm{I}_q}{2\pi}+O\left(\frac{\log(j+1)}{j}\right);\label{c1b}\\
&\E\big(w_n^r(\l_j)\comp{w_{\bs\eps}^s(\l_j)}\big)=\frac{\big(A(\l_j)\big)_{r\bs\cdot}}{2\pi}+O\left(\frac{\log(j+1)}{j\l_j^{d_r^0}}\right), \quad\mbox{ for } j=1,\cdots,m. \label{c1c}
\end{align}
\end{subequations}
\normalsize
By using that $I_T(\l;n)=O(I_n(\l))$, rewrite the expression inside the absolute value on the LHS of \eqref{indir} as
\small
\[ \Big[\big(O(1)w_n^r(\l_{j})-\big(A(\l_j)\big)_{r\bs\cdot}w_{\bs\eps}(\l_j)\big)\Big]\comp{w_n^s(\l_{j})}
+\big(A(\l_j) \big)_{r\bs\cdot}w_{\bs\eps}(\l_j)\Big[\comp{w_n^s(\l_{j})}-\comp{\big(A(\l_j)\big)_{\bs\cdot s}w_{\eps}(\l_j)}\Big].\]
\normalsize
From \eqref{c1a}, \eqref{c1b}, \eqref{c1c}, $\big(A(\l_j)\big)_{r\bs\cdot}\big(\comp{A(\l_j)}\big)_{r\bs\cdot}/2\pi=f^{rr}(\l_j)$ and $f^{rr}(\l_j){\l_j}^{2d_r^0}\sim G^{rr}_0$, it follows that
\small
\begin{align*}
\E\Big(\Big|\big(O(1)w_n^r(\l_{j})&-\big(A(\l_j)\big)_{r\bs\cdot}w_{\bs\eps}(\l_j)\big)\Big|^2\Big)=\\
&=\E\Big(O(1)I_n^{rr}(\l_{j})-O(1)w_n^r(\l_{j}) \comp{w_{\bs\eps}(\l_j)}\big(\comp{A(\l_j)}\big)_{r\bs\cdot}+\big(A(\l_j)\big)_{r\bs\cdot}I_{\bs\eps}(\l_j) \big(\comp{A(\l_j)}\big)_{\bs\cdot s}-\\
&\hskip1.5cm-O(1)\big(\comp{ A(\l_j)}\big)_{r\bs\cdot}w_{\bs\eps}(\l_j)w_n^{r}(\l_{j})\Big)\\
&=O\left(j^{-1}\log\big(j+1\big)\l_{j}^{-2d_r^0}\right)
\end{align*}
\normalsize
and similarly for $\E\big(\big|\comp{w_n^s(\l_j)}-\comp{\big(A(\l_j)\big)_{\bs\cdot s}w_{\eps}(\l_j)}\big|^2\big)$. Also, $\E\big(I_n^{ss}(\l_{j})\big)=O(\l_{j}^{-2d_s^0})$. Finally, the Cauchy-Schwartz's inequality yields
\small\begin{align}
\E(|\mathscr A_j|)&\leq\big|\l_j^{d_r^0+d_s^0}\big| \left[\E\Big(\Big|O(1)w_n^r(\l_j)-\big(A(\l_j)\big)_{r\bs\cdot}w_{\bs\eps}(\l_j)\Big|^2\Big)^{\frac12}\E\big(\big|w_n^s(\l_{j})\big|^2\big)^{\frac12}\ \ +\right.\nonumber\\
&\hspace{1.5cm}+\left.\E\big(\big|\big(A(\l_j) \big)_{r\bs\cdot}w_{\bs\eps}(\l_j)\big|^2\big)^{\frac12}\E\Big(\Big|\comp{w_n^s(\l_j)}-\comp{\big(A(\l_j)\big)_{\bs\cdot s}w_{\bs\eps}(\l_j)}\Big|^2\Big)^{\frac12}\right]\nonumber\\
&=O\left(j^{-1/2}\log(j+1)^{1/2}\right),\nonumber
\end{align}\normalsize
which completes the proof.\fim

\subsubsection{Proof of Corollary \ref{tapered.an}}

By Remark \eqref{rmkan}, it suffices to show that the results in Lemma \ref{lemma1b2}(a) and (b) hold with $f_n(\cdot)=I_T(\cdot;n)$. The results follows by the exactly same argument as in the proof of Lemma 1(b1) and  Lemma 1(b2) of Shimotsu (2007) (see also equations (C2) and (C6)-(C8) in Lobato, 1999) in view of $I_T(\l;n)=O(I_n(\l))$. \fim

\end{document}